\documentclass [smallcondensed] {svjour3}
\usepackage{color}

\usepackage{amsmath,amsfonts}
\usepackage{graphicx}

\newcommand{\ca}[1]{\overline{#1}}
\newcommand{\D}{\mathrm{d}}

\newcommand{\DT}{\mathrm{\Delta}t}

\newcommand{\piu}[1]{{#1}_{j+1/2}^+}
\newcommand{\meno}[1]{{#1}_{j+1/2}^-}

\newcommand{\dx}{\delta}

\begin{document}
\title{Well-balanced high order schemes on non-uniform grids and entropy residuals}
\author{G. Puppo \and M. Semplice}
\institute{
G. Puppo 
\at
Dipartimento di Scienza e Alta Tecnologia
Universit\`a dell'Insubria
Via Valleggio, 11
22100 Como
\email{gabriella.puppo@uninsubria.it}
\and 
M. Semplice
\at
Dipartimento di Matematica ``G. Peano''
Universit\`a di Torino
Via C. Alberto, 10
10123 Torino (Italy)
\email{matteo.semplice@unito.it}
}
\date{Received: date / Accepted: date}
\maketitle

\begin{abstract}
This paper is concerned with the construction of high order schemes on irregular grids for balance laws, including a discussion of an a-posteriori error indicator based on the numerical entropy production. We also impose well-balancing on non uniform grids for the shallow water equations, which can be extended similarly to other cases, obtaining schemes up to fourth order of accuracy with very weak assumptions on the regularity of the grid. Our results show the expected convergence rates, the correct propagation of shocks across grid discontinuities and demonstrate the improved resolution achieved with a locally refined non-uniform grid.

The error indicator based on the numerical entropy production, previously introduced for the case of systems of conservation laws, is extended to balance laws. Its decay rate and its ability to identify discontinuities is illustrated on several tests. The schemes proposed in this work naturally can also be applied to systems of conservation laws.

\end{abstract}

\section{Introduction}
Many problems arising from engineering applications involve the ability to compute flow fields on complex domains, governed by hyperbolic systems of balance laws. Often, many scales are involved and this prompts the need for algorithms that are able to modify the scheme and/or the underlying grid following the evolution of the flow. Several wide purpose codes are available and many of them are based on finite volume schemes, see e.g. Fluent \cite{fluent} or ClawPack \cite{clawpack}. Usually these codes are second order accurate with high order versions, if available, in progress. On the other hand they provide the user with the flexibility of an adaptive grid, which is extremely useful to tackle highly non-homogeneous solutions.

At the same time, high order finite volume schemes are well established in the literature: from the early review in \cite{Shu97} to the more recent paper \cite{Dumbser:DGFV}, extensive studies have been conducted on the construction of high order finite volume schemes. In this paper we carry out a detailed study of the issues arising in finite volume algorithms on irregular grids, and in particular we construct finite volume high order WENO schemes, including the treatment of source terms and addressing the issue of well balancing for steady state solutions. We concentrate on the one-dimensional case, since most problems already arise in this setting. These results can be extended to multidimensional problems discretized with cartesian grids. Schemes based on cartesian grids can be easily parallelized and boundary conditions for complex domains can be implemented with the ghost fluid method as in \cite{Iollo}. 

Adaptive grids can be constructed either by defining a single non uniform grid on which all degrees of freedom are located, as in most unstructured grid managers, or superposing several patches of uniform cartesian grids of different levels of refinement as in the ClawPack solver \cite{clawpack}. In this latter approach the different patches must communicate and the enforcement of conservativity and well balancing for steady states are not straighforward \cite{DonatWellBalanced}. High order schemes for the AMR approach can be found in \cite{BaezaMulet,ShenQiuChristlieb}. For applications to the shallow water equations, see the software GeoClaw \cite{clawpack} and \cite{GeorgeAMRMalpasset}.

In our case we consider a single highly non-uniform grid. Such grids commonly arise in h-adaptive methods \cite{HartenHayman:1983}, expecially when using moving mesh methods \cite{TanqTang:2003,Tang:2004}. In one space dimension, when the grid size varies smoothly, one can remap the problem to a uniform grid as in \cite{FazioLeveque:2003}, but this cannot be expected to work in more space dimensions of when the grid size can jump abruptly as in dyadic/quadtree/octree grid refinement. These latter discretization techniques start from a conforming, often uniform, partitioning of the simulation domain  and allow the local refinement of each control volume by splitting it in $2^d$ parts in $d$ space dimensions, like in \cite{HuGreavesWu:2002:tritree} for simplices and \cite{WangBorthwickTaylor:2004:quadtree} for quads. Lower order schemes on such grids were employed by the authors in \cite{PS:entropy} in one space dimension and in \cite{PS:HYP12} in two space dimensions for general conservation laws. Two-dimensional applications to the shallow water system may be found in \cite{LiangBorthwick:2009:quadtreeSWE}, or in  \cite{Liang}. 

The construction of a fifth order WENO scheme for conservation laws on one-dimensional non-uniform grids, based on the superposition of three parabolas, has been conducted in \cite{WangFengSpiteri}. Here we extend this construction to the case of balance laws, showing how to obtain positive coefficients in the quadrature of the source term. Moreover we also construct a third order scheme based on \cite{LPR01}, characterized by a stencil of three cells. This reconstruction is particularly suited for two-dimensional problems due to its very compact stencil, see \cite{CRS}.

A first key ingredient of this work is the use of semidiscrete schemes which permit to decouple the space from the time discretization: in this fashion the non-uniformity of the grid boils down to an interpolation problem to reconstruct the boundary extrapolated data which interact through the numerical fluxes. Secondly, the use of the Richardson extrapolation as in \cite{NatvigEtAl} is crucial for the preservation of steady states on a non uniform grid, since it allows to enforce equilibrium at the level of each single cell, thus avoiding the need to account for the non-uniformity of the grid. This yields automatic well-balancing over the whole grid, unlike in the block-structured AMR case, where well-balancing has to be enforced not only on each grid patch but also in the projection and interpolation operators that relate the solution on different grid levels \cite{DonatWellBalanced}.

Moreover, we extend the entropy indicator of \cite{PS:entropy} to the case of balance laws. We show that the numerical entropy production provides a measure of the local error on the cell also in the case of balance laws on non-uniform grids.

Before giving the outline of the paper, we briefly introduce the setting and the notation used in the bulk of this work. We consider balance laws with a geometric source term of the form
\begin{equation}\label{e:blaw}
u_t + \nabla \cdot f(u) = g(u,x)
\end{equation}
and we seek the solution on a domain $\Omega$, with given initial conditions.
The computational domain $\Omega$ is an interval, discretized with cells $I_j=(x_{j-1/2}, x_{j+1/2})$, such that $\cup I_j=\Omega$. The amplitude of each cell is $\dx_j= x_{j+1/2}-x_{j-1/2}$, with cell center $x_j =  (x_{j-1/2}+x_{j+1/2})/2$.

We consider semidiscrete finite volume schemes and denote with  $\ca{U}_j(t)$ the cell average of the numerical solution in the cell $I_j$ at time $t$.  The semidiscrete numerical scheme can be written as
\begin{equation}\label{e:semischeme}
\frac{\D}{\D t}\ca{U}_j= - \frac{1}{\delta_j}\left( {F}_{j+1/2}- {F}_{j-1/2}\right) +  G_j(\ca{U},x).
\end{equation}
The numerical fluxes are computed starting from the boundary extrapolated data, namely
\begin{equation}\label{e:fluxes}
{F}_{j+1/2}=\mathcal{F}(U_{j+1/2}^-,U_{j+1/2}^+)
\end{equation}
where $\mathcal{F}$ is a consistent and monotone numerical flux, evaluated on two estimates of the solution at the cell interface $U_{j+1/2}^{\pm}$. These values are obtained with a high order non oscillatory reconstruction, as described in detail in \S \ref{s:reconstruction}.
Finally, $G_j$ is a consistently accurate discretization of the cell average of the source term on the cell $I_j$, see \S \ref{s:wb}.

In order to obtain a fully discrete scheme, we apply a Runge-Kutta method with Butcher's tableau $(A,b)$, obtaining the evolution equation for the cell averages
\begin{equation}
\label{eq:fullydiscrete}
\ca{U}_j^{n+1} = \ca{U}_j^{n} - \frac{\DT}{\delta_j} \sum_{i=1}^s b_i \left(F^{(i)}_{j+1/2}-F^{(i)}_{j-1/2}\right)
+ \DT \sum_{i=1}^s b_i G^{(i)}_j.
\end{equation}
Here $F^{(i)}_{j+1/2}=\mathcal{F}\big(U^{(i),-}_{j+1/2},U^{(i),+}_{j+1/2}\big)$ and the boundary extrapolated data 
$U^{(i),\pm}_{j+1/2}$ are computed from the stage values of the cell averages
\[
\ca{U}_j^{(i)} =
\ca{U}_j^{n} - \frac{\DT}{\delta_j} \sum_{k=1}^{i-1} a_{ik} \left(F^{(k)}_{j+1/2}-F^{(k)}_{j-1/2}\right)
+ \DT \sum_{k=1}^{i-1} a_{ik}G^{(k)}_j.
\]
We point out that the spatial reconstruction procedures of \S \ref{s:reconstruction} and the well-balanced quadratures for the source term of \S \ref{s:wb} must be applied for each stage value of the Runge-Kutta scheme.
In this paper we consider a uniform timestep over the whole grid. A local timestep keeping a fixed CFL number over the grid can be enforced using techniques from \cite{PS:entropy,LambyMullerStiriba}.

We will also consider the preservation of steady state solutions and we will illustrate these techniques on the shallow water system, namely
\begin{equation}\label{eq:swe}
	u=\begin{pmatrix} h\\q\end{pmatrix}
	\qquad
    f(u) = \begin{pmatrix} q\\q^2/h+ \tfrac12 g h^2  \end{pmatrix} 
	\qquad
	g(u,x) = \begin{pmatrix} 0\\-ghz_x\end{pmatrix}
  \end{equation}
Here $h$ denotes the water height, $q$ is the discharge and $z(x)$ the bottom topography, while $g$ is the gravitational constant (see also Figure \ref{fig:sw}). The preservation of steady states depends heavily on the structure of the equilibrium solution one wishes to preserve. Here we will concentrate on the lake at rest solution of the shallow water equation,  given by $H(t,x)=h(t,x)+z(x)=\text{constant}$ and $q(t,x)=0$. Many works have been dedicated to this problem since the paper \cite{BermudezVazquez:1994} shed light on the importance of well-balancing (or C-property). For example, see 
\cite{XingShu:2005:WBSWEfd} in the finite difference setting,
\cite{XingShu:2006:WBDG,NatvigEtAl,NoelleXingShu:2007:SWEmovingwater} in the finite volume setting, 
\cite{XingShu:2006:WBDG,Xing:2013:WBDGmovingwater,CaleffiValiani:2013:RKDG3WB} in the Discontinuous Galerkin framework and
\cite{VignoliTitarevToro:2008:ADERchannel,CastroToroKaser:2012:ADERtsunami} in the ADER setting.

The structure of the paper is as follows: in \S \ref{s:reconstruction} we introduce the third order accurate C-WENO (Compact WENO) reconstruction on non uniform grids, generalizing the results of \cite{LPR01}, and we extend the fifth order accurate WENO reconstruction on non uniform grids of \cite{WangFengSpiteri}, adding the evaluation of the reconstruction at the centre of cells which is needed in the computation of the source term.
In \S \ref{s:wb} we extend the construction of well-balanced schemes of \cite{Audusse:2004,NatvigEtAl} to the non-uniform grid setting.
Next, in \S \ref{s:entropy} we extend the notion of numerical entropy production to non uniform grids for balance laws. Finally, \S \ref{s:numerical} contains numerical tests, which illustrate the consistency between accuracy of the schemes and rate of convergence of the numerical entropy production, for several types of grids.

\section{High order reconstructions on non uniform grids}
\label{s:reconstruction}

The mission of reconstruction algorithms is to give estimates of a function at some points, starting from discrete data. In particular, for finite volume schemes for balance laws, the starting data are the cell averages of a function $v$, and we wish to estimate $v$ at the cell interfaces, and, if needed, at some other internal  points, using a finite dimensional approximation, such as a piecewise polynomial interpolator. Typically, estimates of $v$ at internal points within a cell are needed to compute the cell averages of the source term through a quadrature formula. Thus, the reconstruction will be described as an interpolation algorithm.

Suppose then that we are given the cell averages 
\[
\ca{V}_j = \frac{1}{\dx_j} \int_{I_j} v(x)\; \D x.
\]
of a smooth function $v(x)$. In order to fix ideas, we consider a piecewise polynomial reconstruction 
$\mathcal{R}$ such that
\[
\mathcal{R}(\ca{V},x) = \sum_j \chi_{I_j}(x) P_j(x),
\]
which gives the boundary extrapolated data as
\begin{equation} \label{eq:bdryextrapdata}
V_{j+1/2}^-=P_j(x_{j+1/2}), \qquad V_{j+1/2}^+=P_{j+1}(x_{j+1/2}).
\end{equation}
The reconstruction must be
conservative, i.e.
\[
\frac{1}{\dx_j} \int_{I_j} \mathcal{R}(\ca{V},x)\; \D x = \ca{V}_j,
\]
and high order accurate at the cell interfaces for smooth data, in the sense that
\[
V_{j+1/2}^- = v(x_{j+1/2}) + O(\dx_j)^p, \qquad V_{j-1/2}^+ = v(x_{j-1/2}) + O(\dx_j)^p.
\]
Moreover, the reconstruction should be non-oscillatory, preventing the onset of spurious oscillations. Finally, for accuracy of order higher than 2, the evaluation of the cell average of the source term requires the reconstruction of the point values of $v$ at the  nodes of the well-balanced quadrature formula. For schemes of order 3 and 4, it is enough to reconstruct $v$ at the cell centers, thus we will require that, for smooth $v(x)$, 
\[
V_{j} = v(x_{j}) + O(\dx_j)^p.
\]

\subsubsection*{First order reconstruction}
In this case, the reconstruction is piecewise constant, and we have
\[
V_{j+1/2}^- = \ca{V}_j, \qquad V_{j-1/2}^+ = \ca{V}_j.
\]

\subsubsection*{Second order reconstruction}
Here, the reconstruction is piecewise linear, and we have
\[
V_{j+1/2}^- = \ca{V}_j+\tfrac12 \sigma_j \dx_j, \qquad V_{j-1/2}^+ = \ca{V}_j-\tfrac12 \sigma_j \dx_j,
\]
where $\sigma_j$ is a limited slope, i.e., chosen a limiter
$\Phi$, define the interface slopes as
\begin{equation}\label{e:interface_slope}
\sigma_{j+1/2} = \frac{\ca{V}_{j+1}-\ca{V}_{j}}{x_{j+1}-x_j}=
\frac{\ca{V}_{j+1}-\ca{V}_{j}}{\tfrac12(\dx_j+\dx_{j+1})}  
\end{equation}
then the limited slope within the $I_j$ cell is given by
\[
\sigma_j = \Phi \left( \sigma_{j-1/2},\sigma_{j+1/2}\right).
\]
For a collection of limiting functions, see \cite{LeVeque:book}. In our tests, we have chosen the MinMod limiter.

\subsubsection*{Third order reconstruction}
\begin{figure}
\begin{center}
\includegraphics[width=0.7\textwidth]{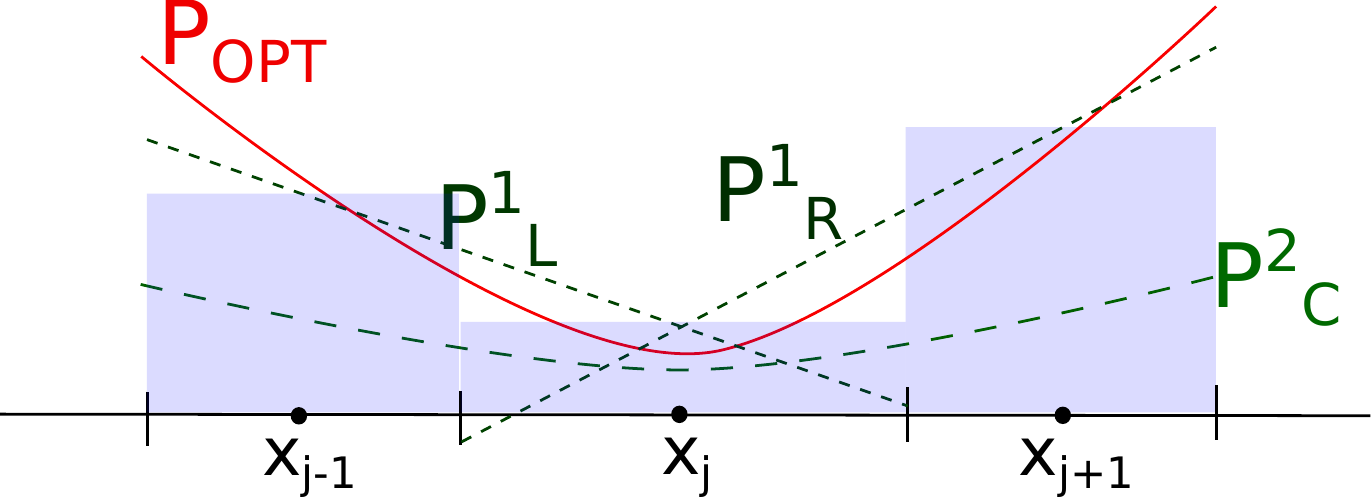}
\end{center}
\caption{\sf Compact WENO reconstruction}\label{f:cweno}
\end{figure}

The third order reconstruction is based on the compact WENO (C-WENO) technique introduced in \cite{LPR01}. This reconstruction is characterized by a particularly compact stencil, which is very important when dealing with adaptive grids. Moreover, unlike the classical WENO third order reconstruction based on the combination of two linear functions, the C-WENO reconstruction contains also a parabola and it remains uniformly third order accurate throughout the interval $I_j$ on smooth flows. To our knowledge, the reconstruction presented here is the first extension of the C-WENO reconstruction to the case of non-uniform grids.
Fig. \ref{f:cweno} illustrates the polynomials composing this reconstruction.

The interpolant is piecewise quadratic, and the  parabola reconstructed in each cell is the convex combination of two linear  functions $P^1_L$, $P^1_R$, and a parabola, $P^2_{C}$. In order to simplify the notation we describe the reconstruction on a reference cell, labelled with the index $j=0$. The two linear functions interpolate $v$ in the sense of cell averages on the stencils $\{I_{-1}, I_0\}$ and $\{I_{0}, I_{+1}\}$. Each of these functions approximates $v$ with order $O(\dx_0)^2$ accuracy uniformly on $I_0$. Further, the parabola $P^2_{\text{OPT}}$ is introduced by the requirement that
\[
\frac{1}{\dx_0}\int_{I_0} P^2_{\text{OPT}}(x) \; \D x = \ca{V}_0, \qquad 
\frac{1}{\dx_{\pm 1}}\int_{I_{\pm 1}} P^2_{\text{OPT}}(x) \; \D x = \ca{V}_{\pm 1}.
\]
This parabola approximates $v$ with order $O(\dx_0)^3$ accuracy uniformly on $I_0$. Next, the parabola $P^2_C$ is introduced, defined as
\[
P^2_{\text{OPT}} = \alpha_0 P^2_C  + \alpha_{+1} P^1_R + \alpha_{-1}P^1_L
\]
with $\alpha_0=\tfrac12$, $\alpha_{\pm 1}=\tfrac14$. 
The reconstruction is given by
\[
P^2(x) = \omega_0 P^2_C  + \omega_{+1} P^1_R + \omega_{-1}P^1_L.
\]
When the function $v$ is smooth, one would like that $\omega_k=\alpha_k + O(\dx_0)^2$, to ensure that $P^2$ has the same accuracy of $P^2_{\text{OPT}}$, otherwise, the non linear weights $\omega_k$ are designed to switch on only the contribution coming from the one-sided stencil on which the function is smooth.

For a non uniform grid, the coefficients of the two linear interpolants on the cell $I_0$ are
\begin{align*}
P^1_R(x)& = \ca{V}_0 + \sigma_{+1/2}(x-x_0) \\
P^1_L(x)& = \ca{V}_0 + \sigma_{-1/2}(x-x_0),
\end{align*}
where $\sigma_{\pm 1/2}$ have been defined in \eqref{e:interface_slope}. The optimal parabola is
\begin{align*}
P^2_{\text{OPT}}  &= a + b(x-x_0) + c(x-x_0)^2, \\
 c&= \frac32 \frac{\sigma_{+1/2}-\sigma_{-1/2}}{\dx_{-1}+\dx_{0}+\dx_{+1}} \\
 b&= \frac{( \dx_{0}+2\dx_{-1})\sigma_{+1/2}+( \dx_{0}+2\dx_{+1})\sigma_{-1/2}}
    {2(\dx_{-1}+\dx_{0}+\dx_{+1})}  \\
 a&= \ca{V}_0  - \tfrac{1}{12}c\, \dx_0^2.
\end{align*}
As in WENO-like reconstructions, the non linear weights $\omega_k$ are computed
as
\[
\tilde{\omega}_k = \frac{\alpha_k}{(\epsilon + \text{IS}_k)^2}, \qquad \omega_k=\frac{\tilde{\omega}_k}{\sum_{l=-1}^1 \tilde{\omega}_l},
\]
starting from the smoothness indicators $\text{IS}_k$ defined in \cite{Shu97}. In this case, they are given by
\begin{align*}
\text{IS}_{-1} & = \dx_0^2 \sigma_{-1/2}^2 \\
\text{IS}_1 & = \dx_0^2 \sigma_{+1/2}^2 \\
\text{IS}_0 & = \frac{1}{\alpha_0^2} \left[  
\left( b-\alpha_{-1}\sigma_{-1/2}-\alpha_{+1}\sigma_{+1/2}\right)\dx_0^2
+ \tfrac{13}{3}c^2\, \dx_0^4
\right].
\end{align*}
Since $P^2_{\text{OPT}}$ is uniformly third order accurate on the whole interval, the boundary extrapolated data and the value $V_0$ at the cell center are all computed evaluating the same quadratic polynomial at the corresponding points inside the cell.

\subsubsection*{Fourth order reconstruction}
\begin{figure}
\begin{center}
\includegraphics[width=0.7\textwidth]{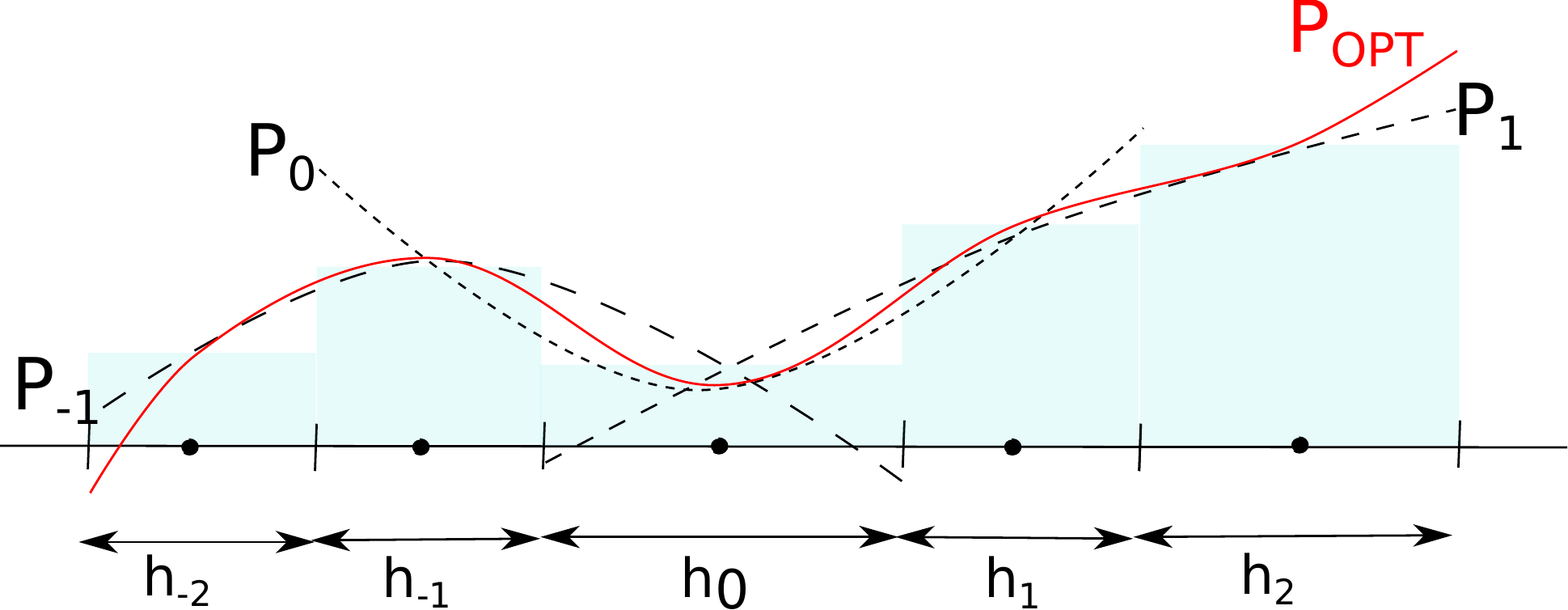}
\end{center}
\caption{\sf Parabolic WENO reconstruction}\label{f:WENO5}
\end{figure}

The fourth order reconstruction is based on the fifth order WENO reconstruction computed from the convex combination of three parabolas, as in \cite{Shu97}. The coefficients of the combination of the three parabolas are computed in order to yield fifth order accuracy at the boundary of the cell, see Fig \ref{f:WENO5}. It is tedious but straightforward to see that positive coefficients can be found to result in fifth order accuracy at the cell interfaces even on non uniform grids (see below and \cite{WangFengSpiteri}). However, there is no set of positive coefficients resulting in fifth order accuracy at the cell center, see \cite{NatvigEtAl}. Here we show that it is possible to find three positive coefficients giving {\em fourth} order accuracy at the center of the cell.

For the sake of completeness, we review the coefficients of the reconstruction on non uniform grids, as in \cite{WangFengSpiteri}, using the notation established in Fig. \ref{f:WENO5}. Again we consider a reference cell with index $0$. The goal of the reconstruction is to mimic the quartic polynomial $P_{\text{OPT}}$ interpolating the data $\ca{V}_{l}, l=-2, \dots, 2$ in the sense of cell averages. Clealy,  $P_{\text{OPT}}$ would provide fifth order accuracy uniformly in the interval $I_0$, in the case of smooth data.

For each point $\hat{x}$ in which the reconstruction is needed, we look for three positive coefficients $d_{-1}, d_{0}, d_{1}$ that add up to $1$ and  such that 
\begin{equation}\label{eq:WENO5:def}
P_{\text{OPT}} (\hat{x}) = \sum_{l=-1}^1 d_l P_l(\hat{x}),
\end{equation}
where the $P_l$'s are the three parabolas, interpolating in the sense of cell averages the data $\ca{V}_{l-1}, \ca{V}_{l}, \ca{V}_{l+1}$. The coefficients of the three parabolas can be found in \cite{WangFengSpiteri}. Here we give the linear weights that permit to reconstruct the left and right boundary extrapolated data. To simplify the notation, we write
\begin{equation} \label{eq:notazionebrutta}
\dx_l^k = \sum_{i=l}^k \dx_i, 
\end{equation}
then the coefficients for the boundary extrapolated data $V_{+1/2}^-$ are
\begin{align*}
d_{1} &= \frac{\dx_{-1}(\dx_{-2}+\dx_{-1})}{\dx_{-2}^2 \dx_{-1}^2} 
\\
d_0 &= \frac{\dx_{0}^2(\dx_{-2}+\dx_{-1})(\dx_{-2}^1 + \dx_{-1}^2)}
{\dx_{-2}^2 \dx_{-1}^2 \dx_{-2}^1} 
\\
d_{-1} &= 
\frac{\dx_{0}^2(\dx_{0}+\dx_1)}{\dx_{-2}^2 \dx_{-2}^1} 
\end{align*}
Note that, if  $\dx_{-2}=\dx_{-1}=\dx_0=\dx_1=\dx_2$, then  
$d_{-1}= \tfrac{3}{10},  d_0= \tfrac{3}{5}, d_1= \tfrac{1}{10}$, as in the usual uniform grid case. Similarly, the coefficients for the reconstruction of $V_{-1/2}^+$ are
\begin{align*}
d_{-1} &= \frac{\dx_1(\dx_1+\dx_2)}{\dx_{-2}^2 \dx_{-2}^1} 
\\
d_0 &= \frac{\dx_{-2}^0(\dx_1+\dx_2)(\dx_{-2}^1 + \dx_{-1}^2)}
{\dx_{-2}^2 \dx_{-1}^2 \dx_{-2}^1} 
\\
d_1 &= 
\frac{\dx_{-2}^0(\dx_{-1}+\dx_0)}{\dx_{-2}^2 \dx_{-1}^2} 
\end{align*}
We remark that the coefficients $d_k$ are positive and add up to $1$, so that \eqref{eq:WENO5:def} is a convex combination, for all possible values of the local grid size $\dx_{-2},\ldots,\dx_{2}$.

For the $5^{\text{th}}$-order reconstruction at cell center $x_0$, one 
finds negative coefficients even for uniform meshes. In fact, see \cite{NatvigEtAl},
$d_{-1} =-\tfrac{9}{80}, d_0 =\tfrac{49}{40}, d_1 =-\tfrac{9}{80}$. Since the well balanced quadrature based on the three points $x_{\pm 1/2}, x_0$ is only fourth order accurate, there is actually no need for fifth order accuracy in this case. Thus, we look for {\em positive} coefficients $d_0, d_{\pm 1}$ such that $1=\sum d_l$, and $V_0$ is fourth order accurate,
\[
V_0 = \sum_{l=-1}^1 d_l P_l(x_0) = v(x_0) + O(\dx_0)^4.
\]
After tedious computations, we find that $d_1$ and $d_{-1}$ must satisfy
\[
\dx_{-2}^1 d_{-1} - \dx_{-1}^2d_1 = \dx_1 - \dx_{-1}
\]

\begin{figure}
\begin{center}
\includegraphics[width=0.3\textwidth]{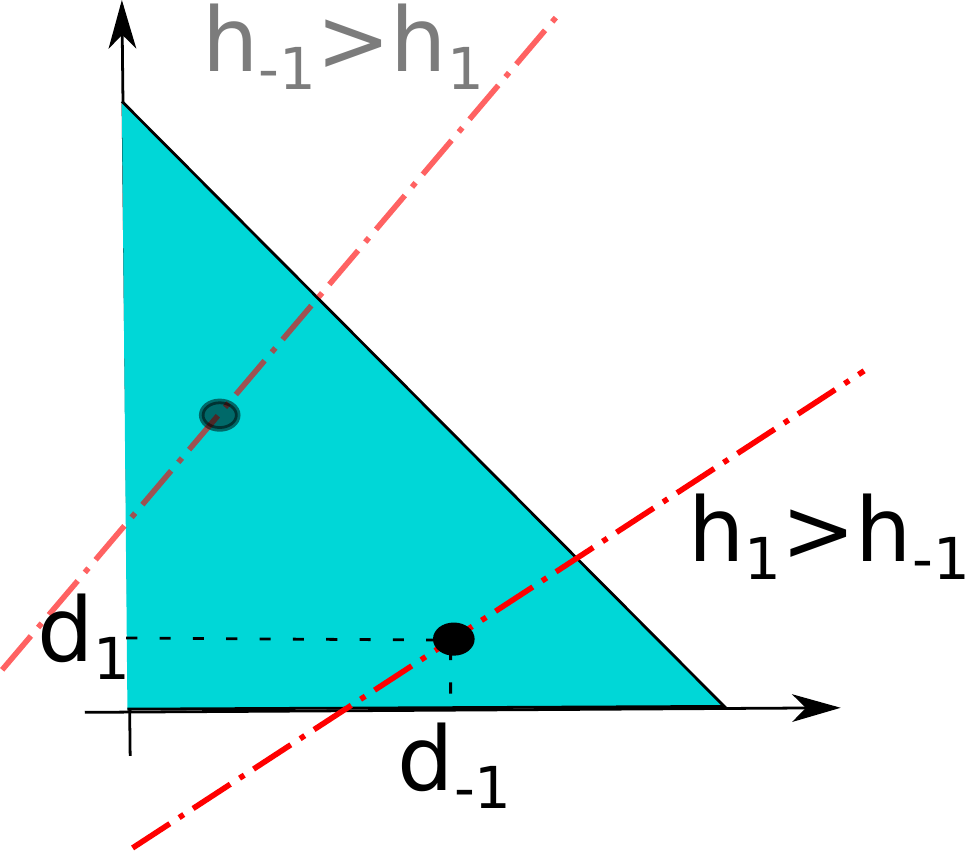}
\end{center}
\caption{\sf Reconstruction of the point value in the cell center for  WENO. Locus of positive linear weights (dash-dot lines) and the coefficients chosen by \eqref{eq:CWEN04:center} (black dots).}\label{fig:WENO4}
\end{figure}
\noindent Since we wish all coefficients to be positive, the solution must be sought in the simplex shown in Fig. \ref{fig:WENO4}. Clearly, the solution is over-determined, we pick the values that maximize the size of the minimum coefficient, that is
\begin{equation}\label{eq:CWEN04:center}
\text{If } \dx_1 > \dx_{-1} \left\{
\begin{aligned}
& d_1 = \frac12
\frac{\dx_{-2} + 2 \dx_{-1}+\dx_0}{\dx_{-2}^1 + \dx_{-1}^2}  \\
& d_{-1} = \frac{\dx_1 - \dx_{-1} + d_1 \dx_{-1}^2}{\dx_{-2}^1}
\\
& d_0 = 1 - d_{-1} - d_1
\end{aligned} \right. ,
\qquad 
\text{else } \left\{
\begin{aligned}
& d_{-1} = \frac12
\frac{\dx_{2} + 2 \dx_{1}+\dx_0}{\dx_{-2}^1 + \dx_{-1}^2}  \\
& d_{1} = \frac{\dx_{-1} - \dx_{1} + d_{-1} \dx_{-2}^1}{\dx_{-1}^2}
\\
& d_0 = 1 - d_{-1} - d_1
\end{aligned} \right.
\end{equation}
where again we have used the convention \eqref{eq:notazionebrutta}.

\section{Well-balanced schemes}
\label{s:wb}
It is important to perform numerical integration of a system of balance laws with schemes that preserve the steady states exactly at a discrete level (well-balanceed schemes), since only these allow to distinguish small perturbations of these states from numerical noise \cite{BermudezVazquez:1994}.

In this section we describe a technique to obtain well-balanced schemes on non-uniform grids for the shallow water equations, with particular attention to the lake at rest solution. In this case, beside well-balancing, it is also particularly important to preserve the positivity of the water height.
We use and generalize to nonuniform meshes the techniques of \cite{Audusse:2004} for obtaining well-balanced schemes irrespectively of the chosen numerical fluxes and of \cite{NatvigEtAl} to obtain high order accuracy through Richardson extrapolation.

\begin{figure}
\begin{center}
\includegraphics[width=0.7\linewidth]{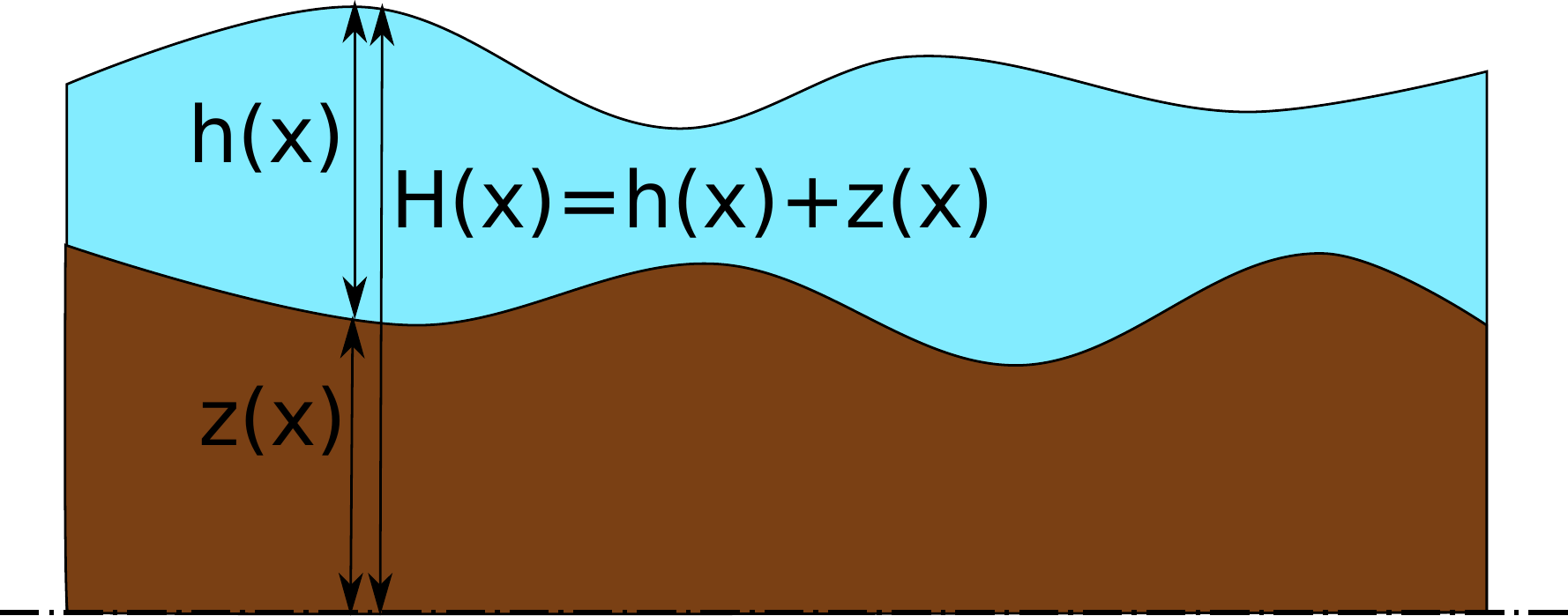}
\end{center}
\caption{Shallow water set up.}
\label{fig:sw}
\end{figure}

There are two sources of error in well-balanced schemes. We illustrate them with a very simple example. We consider a first order reconstruction with the Lax-Friedrichs numerical flux on the lake at rest solution (see Fig. \ref{fig:sw} for notation), thus we suppose that for every index $j$, $q^n_j=0$ and $h^n_j+z_j=H$. The discretized equation on a uniform grid would be
\begin{align*} 
h_j^{n+1} &= h_j^n + \tfrac{\lambda}2\alpha \left(h^n_{j+1} - 2h^n_j +h^n_{j-1}\right)\\
q_j^{n+1} &= -\tfrac{\lambda}4 g\left( (h^n_{j+1})^2 - (h^n_{j-1})^2 \right)
+\tfrac{\lambda}2 gh^n_j \left( z_{j+1}- z_{j-1}\right)
\end{align*}
where we have already substituted $q^n_j=0$. It is easy to see that in the first equation, $h$ does not remain constant because the artificial diffusion term introduces a perturbation whenever $z(x)$ is not constant. In order to prevent this kind of perturbation it is enough to reconstruct along equilibrium variables or to ensure that the boundary extrapolated values at the interface are continuous when equilibrium occours. In the second equation, the perturbation due to the artificial diffusion does not appear exactly because $q$ is an equilibrium variable for the lake at rest equilibrium. However there is a lack of balance betweeen the source and the fluxes at the discrete level: in fact one finds that $q_j^{n+1}=-\tfrac{\lambda}{4}(z_{j+1}^2-2z_jz_{j+1}+2z_jz_{j-1}-z_{j-1}^2)$, which is in general nonzero, unless the bottom is flat.

For these reasons we use the hydrostatic reconstruction of \cite{Audusse:2004} which ensures that the reconstruction is continuous across interfaces when the system is in equilibrium and moreover preserves positivity of the water height. Given a reconstruction algorithm $\mathcal{R}$ with accuracy of order $p$, reconstruct the equilibrium variables $H$ and $q$, obtaining the boundary extrapolated data as in equation \eqref{eq:bdryextrapdata}. In order to ensure that the water height appearing in the fluxes remains non-negative, one locally modifies the bottom by computing boundary extrapolated data also for $h$ and defining
\[
z_{j+1/2}^{\pm} = H_{j+1/2}^{\pm} - h_{j+1/2}^{\pm}
\]
and these are used to compute the bottom topography at the interface
\[
z_{j+1/2} = \max(\piu{z},\meno{z}).
\]
Once these are known, the interface values of $h$ are corrected giving new values
\[
\widehat{h}_{j+1/2}^{\pm} = \max(H_{j+1/2}^{\pm} - z_{j+1/2},0).
\]
Note that $\widehat{h}_{j+1/2}^{\pm}\geq0$ and that at equilibrium $\widehat{h}_{j+1/2}^{+}=\widehat{h}_{j+1/2}^{-}$. The numerical fluxes \eqref{e:fluxes} are then applied to the states 
\[ U^{\pm}_{j+1/2} = \left[\widehat{h}_{j+1/2}^{\pm}, \;
\widehat{h}_{j+1/2}^{\pm} v_{j+1/2}^{\pm} \right]. \]
Here $v_{j+1/2}^{\pm}$ denotes the velocity, obtained as $v_{j+1/2}^{\pm}=q_{j+1/2}^{\pm}/\widehat{h}_{j+1/2}^{\pm}$ or through a desingularization procedure as proposed in \cite{Kur:desing}. Since the reconstruction is continous at equilibrium, for lake at rest data, for each consistent numerical flux, one has ${\mathcal F}(U^-_{j+1/2},U^+_{j+1/2})=f(U^{\pm}_{j+1/2})$. In this fashion Audusse et al.  are able to ensure well-balancing independently on the particular numerical flux used \cite{Audusse:2004}.

In order to complete the semidiscrete scheme \eqref{e:semischeme} we still need to specify the discretization of the source term. For a first order scheme it is enough to choose
\begin{equation}\label{eq:S:1}
G_j = \frac{g}2 \begin{pmatrix}
0\\
(\widehat{h}^-_{j+1/2})^2 - (\widehat{h}^+_{j-1/2})^2
\end{pmatrix}.
\end{equation}
Note that at equilibrium, the above expression exactly cancels out the numerical fluxes and thus the lake at rest solution is preserved at the discrete level. Consistency is obtained through the dependence of $\widehat{h}$ on $z$.

At second order, the second component of the source term is
\begin{align}\label{eq:S:2}
G_{j,2} = \frac{g}2 
 & (\, (\widehat{h}^-_{j+1/2})^2 -(h^-_{j+1/2})^2 
 + (h^+_{j-1/2}+h^-_{j+1/2})(z^+_{j-1/2}-z^-_{j+1/2}) \\
 &  + (h^+_{j-1/2})^2  -(\widehat{h}^+_{j-1/2})^2\, ) \nonumber
\end{align}

On the lake at rest solution, the two $\widehat{h}$ terms cancel the numerical fluxes, while the other terms add up to zero, again giving a well-balanced scheme \cite{Audusse:2004}. On the other hand, off equilibrium, the first and the last two terms cancel by consistency and the middle term is consistent with the cell average of the source. Clearly, equation \eqref{eq:S:2} must be applied to both of the stages of the second order Runge-Kutta method needed to achieve second order accuracy also in time.

For higher orders, we use Richardson extrapolation as in \cite{NatvigEtAl}. This technique is particularly useful on non-uniform grids because it concentrates all the computational effort for the source term within one cell. In fact, the subcell resolution required to compute the quadrature of the source term with high order accuracy can be naturally applied introducing uniformly distributed nodes within each cell. Thus the high order evaluation of the source term is performed entirely within one cell and the coefficients of the quadrature formula will not be affected by the nonuniformity of the mesh.
The source can be rewritten as
\begin{equation}\label{eq:S:4}
G_j = \frac{g}2 \begin{pmatrix}
0\\
(\widehat{h}^-_{j+1/2})^2 -(h^-_{j+1/2})^2 
 + \widetilde{G}_j 
 + (h^+_{j-1/2})^2  -(\widehat{h}^+_{j-1/2})^2 
\end{pmatrix}.
\end{equation}
At second order, 
\[ \widetilde{G}_j = (h^+_{j-1/2}+h^-_{j+1/2})(z^+_{j-1/2}-z^-_{j+1/2}) 
       = \int_{x_{j-1/2}}^{x_{j+1/2}} hz_x \mathrm{d}x + O(\dx_j^2).\]
For order up to four, it is enough to choose
\begin{align*}
 \widetilde{G}_j = & \frac43 \left( 
      (h^+_{j-1/2}+h_{j})(z^+_{j-1/2}-z_j) 
+     (h_j+h^-_{j+1/2})(z_j-z^-_{j+1/2}) 
 \right) \\
& -\frac13 (h^+_{j-1/2}+h^-_{j+1/2})(z^+_{j-1/2}-z^-_{j+1/2}) ,
\end{align*}
where $h_j$ and $z_j$ denote the reconstruction at the center of the cell, which is why we have developed high order reconstructions for the point values of the solution in $x_j$.
Again, equation \eqref{eq:S:4} will be applied to all stages of the Runge-Kutta method used in the fully discrete scheme.


\section{Numerical entropy production for balance laws}
\label{s:entropy}
We wish to devise an error indicator for driving adaptive schemes for balance laws. In particular we extend the notion of numerical entropy production proposed in \cite{P:entropy,PS:entropy} to the case of balance laws with a geometric source term.

In the homogeneous case, that is for systems of hyperbolic conservation laws, the entropy is defined as a convex function $\eta(u)$ for which there exists a function $\psi(u)$ (called entropy flux) such that $\nabla^T\eta f' = \nabla^T\psi$ where $f'$ denotes the Jacobian of the flux function $f$. Then, on smooth solutions, 
\[ \partial_t\eta+\partial_x\psi=0,\]
while on entropic shocks 
\[\partial_t\eta+\partial_x\psi\leq 0\]
in a weak sense, thus singling out the correct unique solutions \cite{Dafermos}.
One can exploit this structure at the discrete level to devise a regularity indicator for finite volume schemes for conservation laws. A fully discrete finite volume conservative scheme for a hyperbolic system can be written in the form 
\[
\ca{U}^{n+1}_j= \ca{U}^{n}_j - \lambda \left( F_{j+1/2}- F_{j-1/2}\right).
\]
Here 
\[
F_{j+1/2} = \sum_{i=1}^s b_i \mathcal{F}\left(U^{(i),-}_{j+1/2},U^{(i),+}_{j+1/2}\right),
\]
$\mathcal{F}$ is a consistent and monotone numerical flux and $U^{(i),\pm}_{j+1/2}$ denote the boundary extrapolated data computed on the $i$-th stage value.

Choosing a numerical entropy flux $\mathcal{P}$, consistent with the exact entropy flux $\psi$, we can define the quantity
\begin{equation} \label{eq:S}
S^n_j = \frac{1}{\DT_n} \left[
\ca{\eta(U^{n+1})}_j - \ca{\eta(U^{n})}_j
+\lambda \left(P_{j+1/2}-P_{j-1/2}\right)
 \right]
\end{equation}
where
\[
P_{j+1/2} = \sum_{i=1}^s b_i \mathcal{P}\left(U^{(i),-}_{j+1/2},U^{(i),+}_{j+1/2}\right)
\]

In \cite{PS:entropy} we proved that 
\[
S^n_j=
\begin{cases}
O(h^p) & \text{on smooth flows} \\
\sim C/h & \text{on shocks}
\end{cases}
\]
where $C$ does not depend on $h$ and $p$ is the order of accuracy of the scheme.
Moreover, if the numerical flux can be written in viscous form as
\[
\mathcal{F}(U^-,U^+) = \tfrac12 (f(U^-)+f(U^+)) - \tfrac12 Q(U^-,U^+)\, (U^+-U^-)
\]
we choose the numerical entropy flux as
\begin{equation}\label{eq:numentflux}
\mathcal{P}(U^-,U^+) = \tfrac12 (\psi(U^-)+\psi(U^+))- \tfrac12 Q(U^-,U^+)\, (\eta(U^+)-\eta(U^-)).
\end{equation}
Then we see numerically that the numerical entropy production is essentially negative definite on smooth flows, in the sense that positive values of $S_j^n$ may occour near local extrema, but their amplitude decreases faster than the order of convergence of the scheme. In particular, we have proved this claim for the upwind and Lax Friedrichs numerical flux applied to first order schemes in the scalar case \cite{PS:entropy}.

We wish to extend this construction to systems of $n$ balance laws. In the case of separable balance laws in the sense of \cite{XingShu:2006:WBDG}, namely if the source can be written as
\begin{equation} \label{eq:separable}
g(u,x) = \sum_{j=1}^M s_j(u,x)z'_j(x)
\end{equation}
(with $s_j:\mathbb{R}^n\times\mathbb{R}\to\mathbb{R}^n$),
the balance law can be rewritten as an homogeneous system of $n+M$ equations. For the case $M=1$, denoting with $A(u)$ the $n\times n$ Jacobian matrix of the flux $f$, one has
\begin{equation}
\label{eq:M1}
\partial_t \begin{pmatrix}u\\z_1\end{pmatrix}
+
\begin{pmatrix}
A(u) & s_1(u,x)\\0&0
\end{pmatrix}
\partial_x \begin{pmatrix}u\\z_1\end{pmatrix}
=
0.
\end{equation}
Exploiting this structure one can extend the notion of entropy. In fact the entropy-entropy flux pair for the balance law must satisfy 
\begin{equation}\label{e:entropy_fluxes}
\left[
\nabla^T_u\eta A(u) ,\; \nabla^T_u\eta  \cdot s_1(u,x) 
\right]
=
\left[
\nabla^T_u\psi ,\; \partial_{z_1} \psi 
\right]
\end{equation}
Note that the $z$-derivative of $\eta$ does not appear in the compatibility condition above, and thus convexity with respect to $z$ is not required. This construction can be easily extended for $M>1$.

Thus we still have entropy conservation for the balance law in the smooth case, provided the entropy-entropy flux pair satisfies \eqref{e:entropy_fluxes}, and the entropy residual defined in \eqref{eq:S} gives a measure of the local error of the numerical scheme.

In the shallow water case, the entropy pair can be chosen as
\begin{equation}\label{eq:shentropy}
\eta(h,u) = \tfrac12 \left(hu^2+gh^2\right) +ghz
\qquad
\psi(h,u) = \eta(h,u)u + \tfrac12 gh^2u,
\end{equation}
see \cite{Bouchut:book}. Note that the function $\eta$ represents the total energy of the system including the potential energy due to the bottom topography. In the following section we will show that the entropy residual converges with the expected rate on smooth flows and detects the presence of shocks in the solution.

\section{Numerical tests}
\label{s:numerical}
The following tests asses the accuracy of the high order reconstructions on non-uniform grids proposed in this work, the well-balancing properties of the fully discrete schemes for the shallow water equations, the resolution of discontinuities on non-uniform grids and the performance of the entropy residual as an error indicator.

In all tests we used the local Lax-Friedrichs numerical flux and the entropy residual defined with the corresponding numerical entropy flux  \eqref{eq:numentflux}, unless otherwise stated.

\begin{figure}
\includegraphics[width=0.45\linewidth]{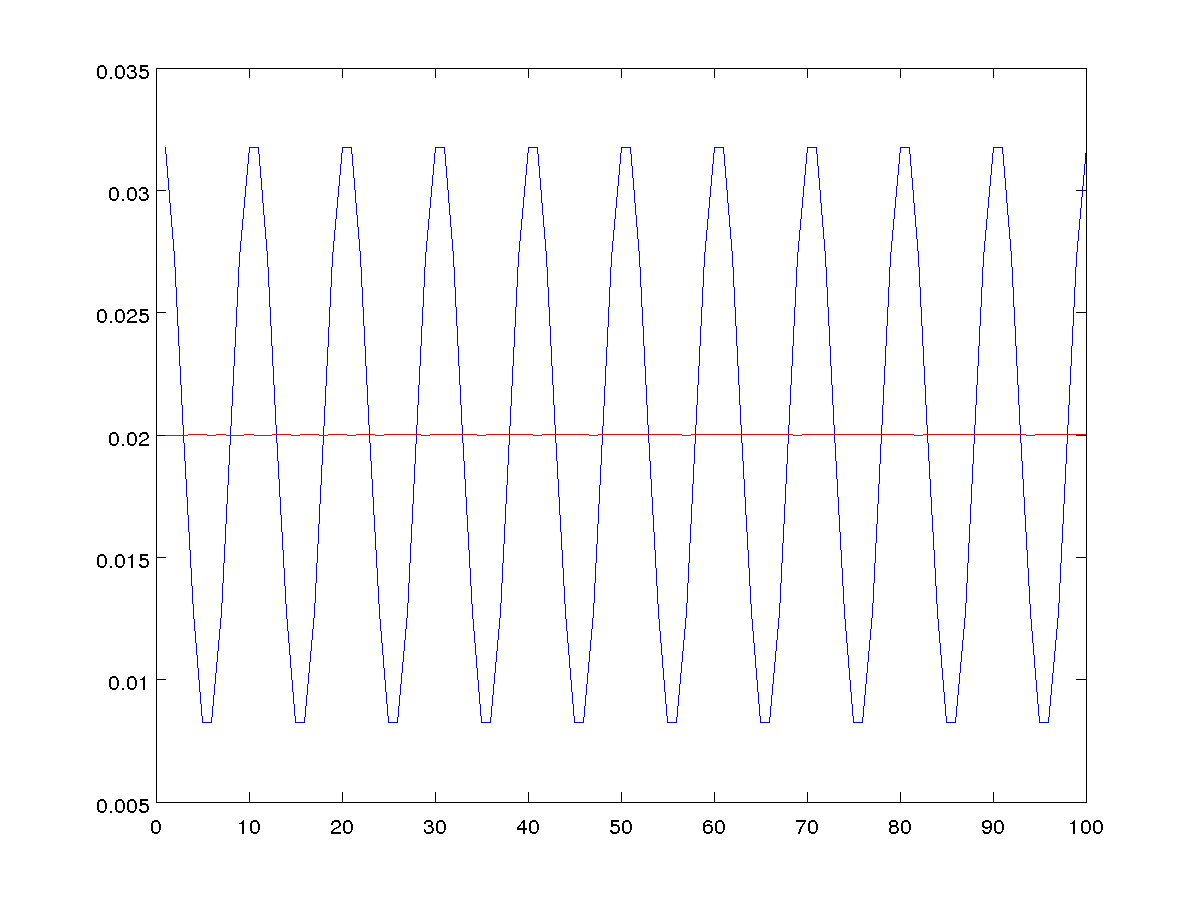}
\hfill
\includegraphics[width=0.45\linewidth]{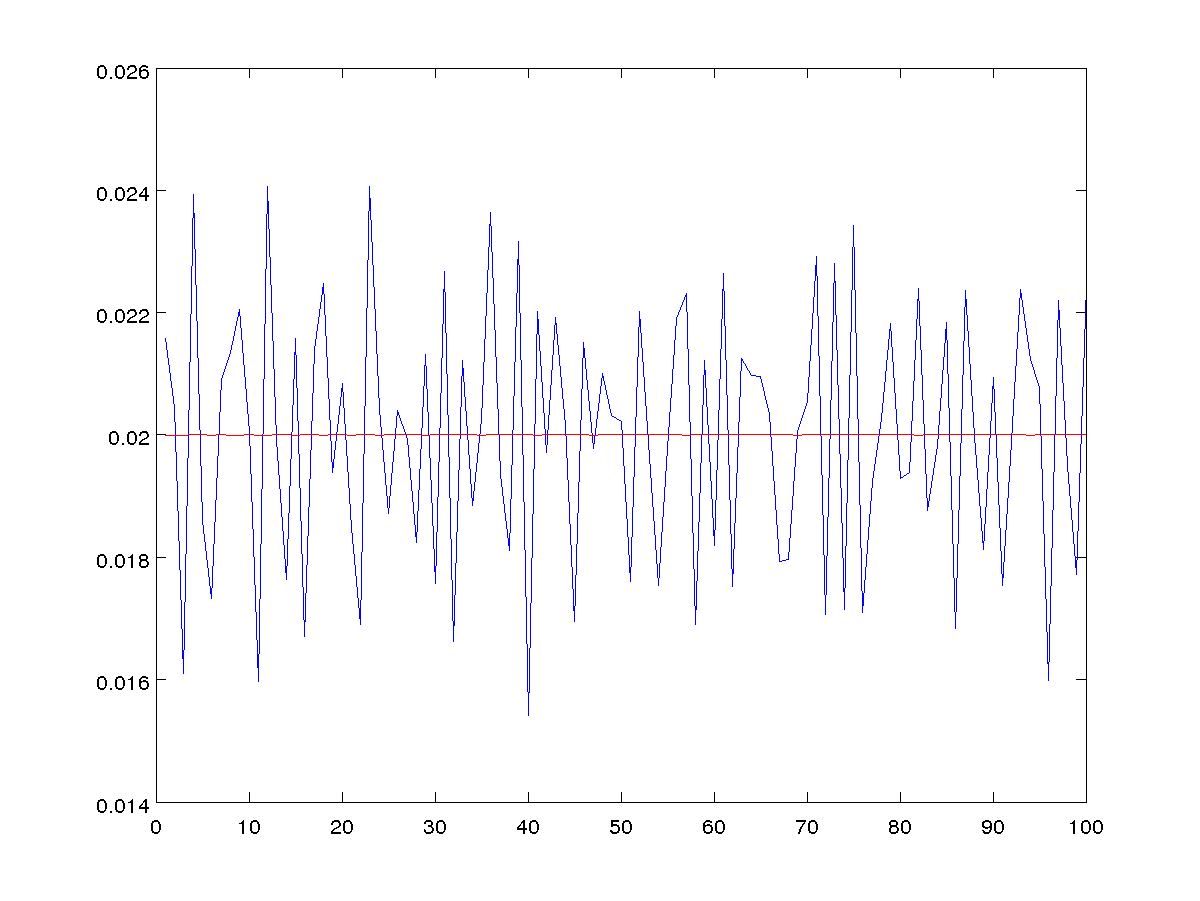}
\caption{Grid spacing for the nonuniform grids used in the numerical tests, shown for the case of $100$ points in $[0,2]$. {\em Quasi-regular} grids (left) and {\em random} grids (right). 
}
\label{fig:grids}
\end{figure}

\paragraph{Grids}
In the numerical tests we use several grids that will be referred to as {\sl uniform}, {\sl quasi-regular}, {\sl random} and {\sl locally refined}. 
For simplicity we define them on the reference interval $[0,1]$.
The quasi regular grid is obtained as  the image of a uniform grid with spacing $\delta=1/N$ under the map
\[\varphi(x)=x+0.1*\sin(10\pi x)/5;\]
The resulting grid spacing is depicted in the left panel of Figure \ref{fig:grids}: we point out that 
\[ 
 (1-\tfrac{\pi}{5})\tfrac1N \leq \dx_j \leq  (1+\tfrac{\pi}{5})\tfrac1N.
\]

Next, we consider non-uniform rough grids that are obtained moving randomly the interfaces of a uniform grid, namely starting from a uniform grid with spacing $\delta$ we consider grids with interfaces at
\[ \tilde{x}_{j+1/2} = j\dx+ \xi_j\tfrac{\dx}{4}\]
where $\xi_j$ are random numbers uniformly distributed in $[-0.5,0.5]$. A realization of such a grid is shown in 
the right panel of Figure \ref{fig:grids}. Here it is easily seen that
\[ 
\tfrac34 \tfrac1N \leq \dx_j \leq  \tfrac54 \tfrac1N.
\]
We use this grid for the purpose of illustration even if of course one would not use such an irregular grid in an application. This grid will be referred to as {\em random} grid.

In some tests we need a grid which is locally refined around a given point $w_C$. For this purpose we consider a grid which, on the standard domain $[0,1]$ is a map of a uniform grid under the function
\begin{equation} \label{eq:locallyrefined}
\varphi(w) = w + 3w(1-w)(w_C-w);
\end{equation}
where $w_C$ is the location in $[0,1]$ of the point where the grid should have its minimum spacing (see e.g. Fig. \ref{fig:steady:transshock:3}).

\subsection{High order schemes on non-uniform grids}
\begin{figure}
\begin{tabular}{cc}
First order
&
Second order
\\
\includegraphics[width=0.45\linewidth]{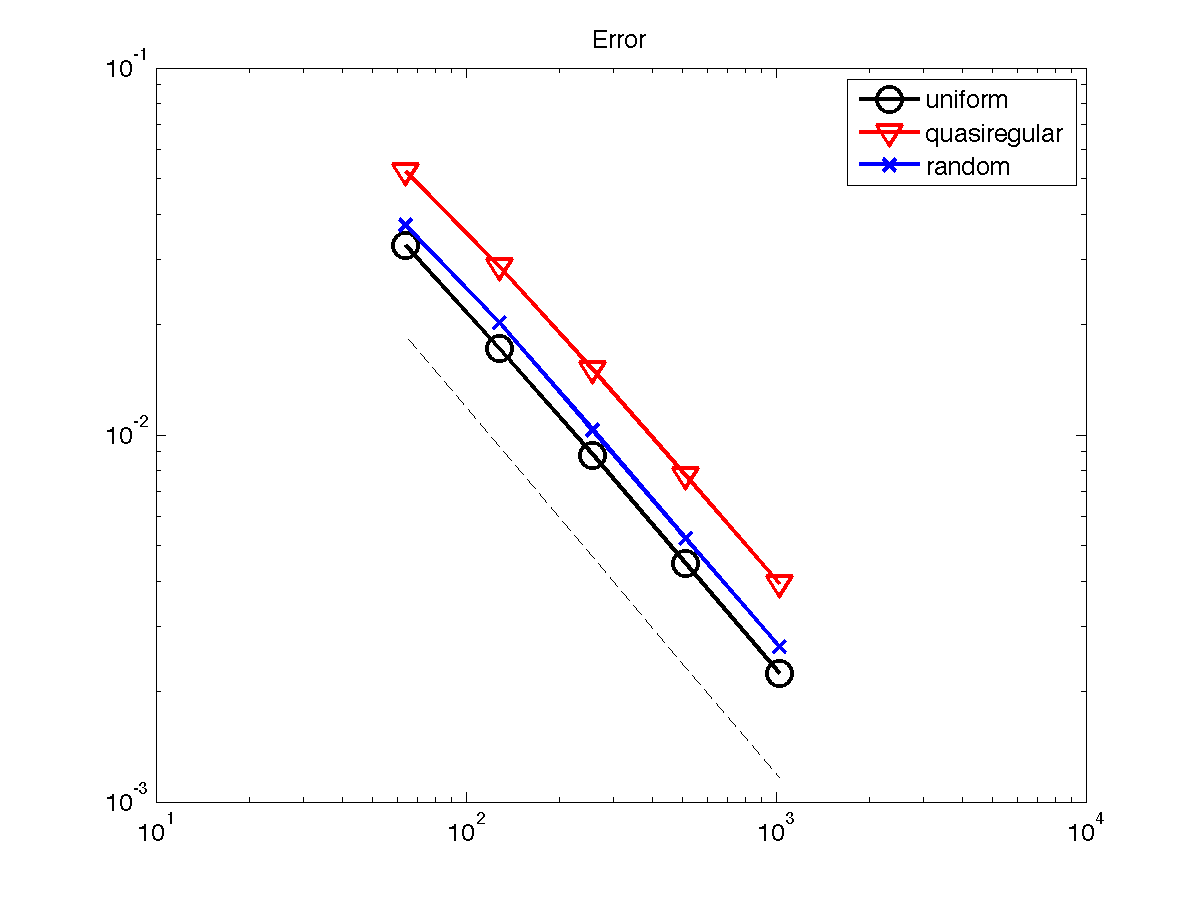}
&
\includegraphics[width=0.45\linewidth]{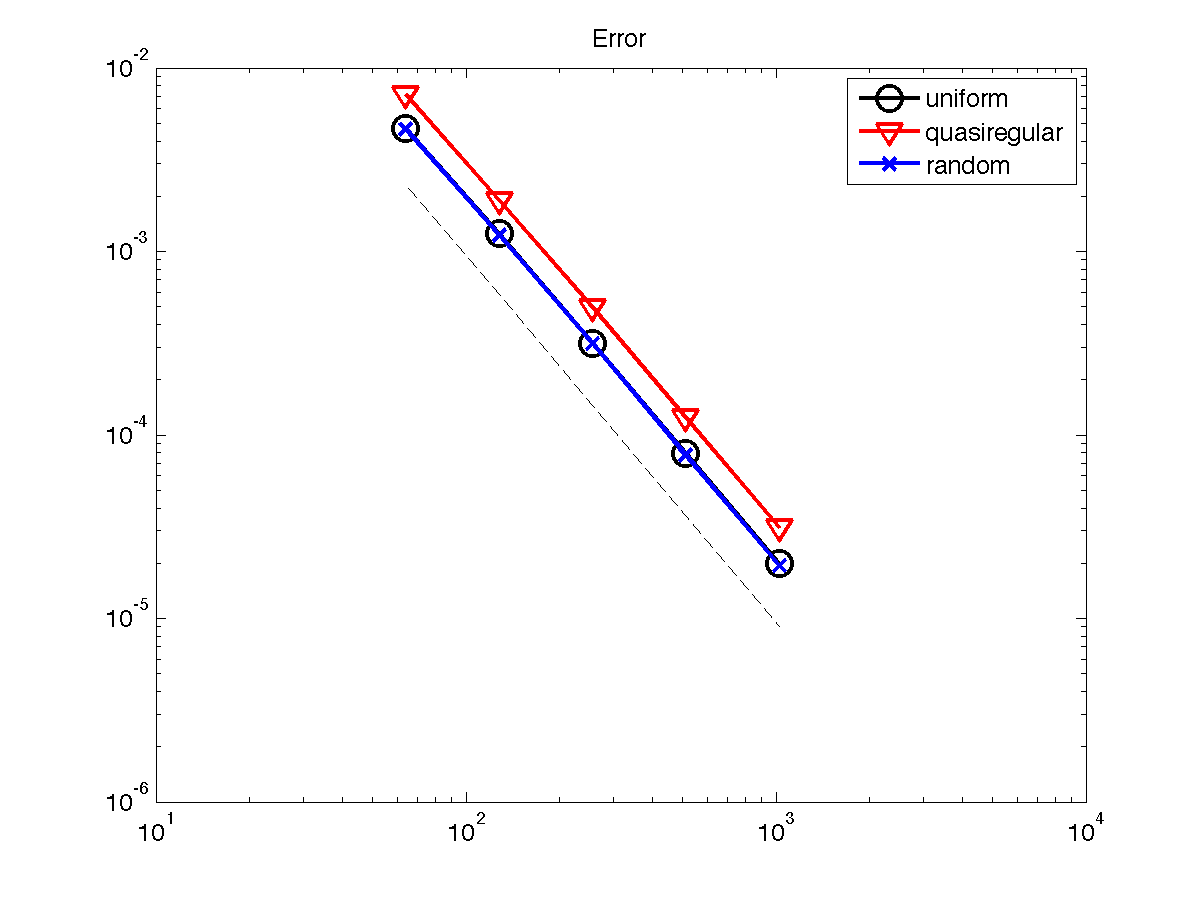}
\\
Third order
&
Fourth order
\\
\includegraphics[width=0.45\linewidth]{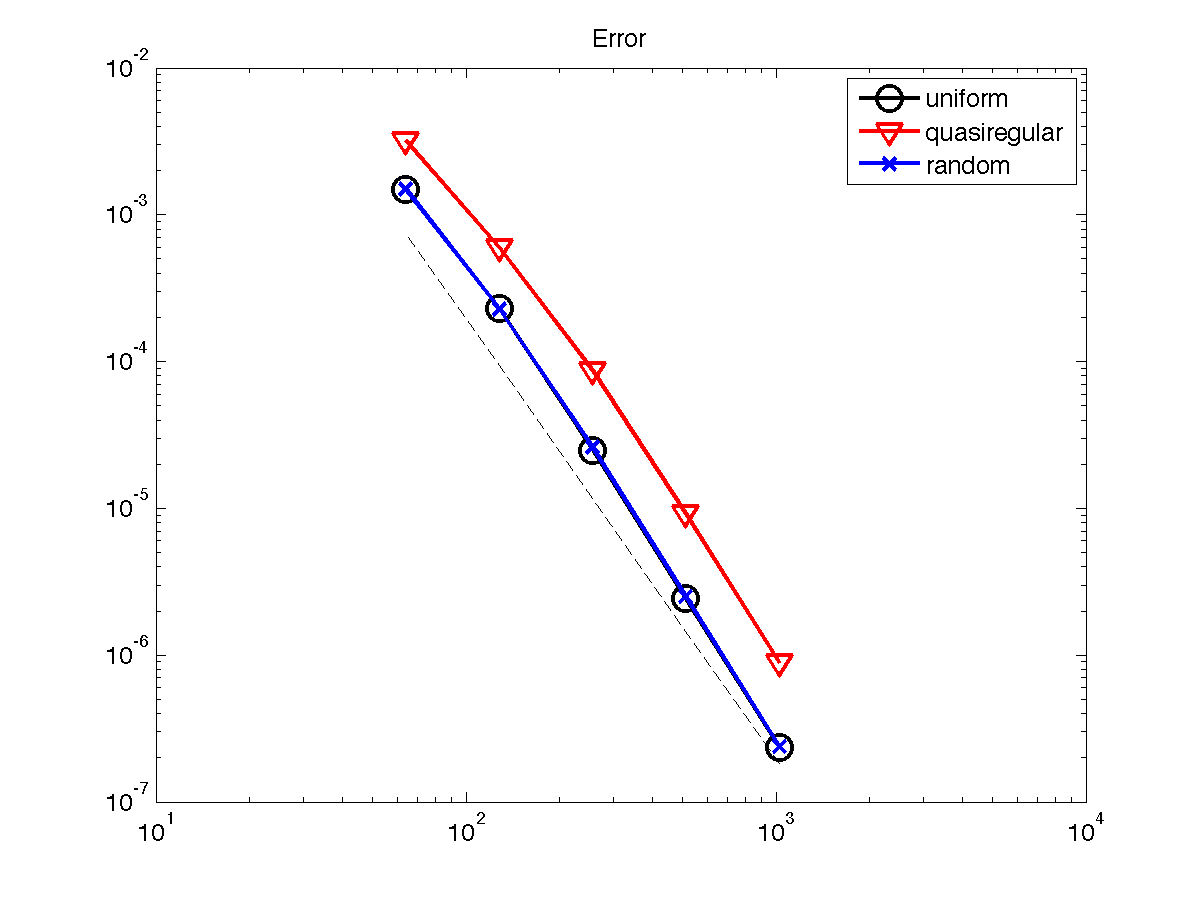}
&
\includegraphics[width=0.45\linewidth]{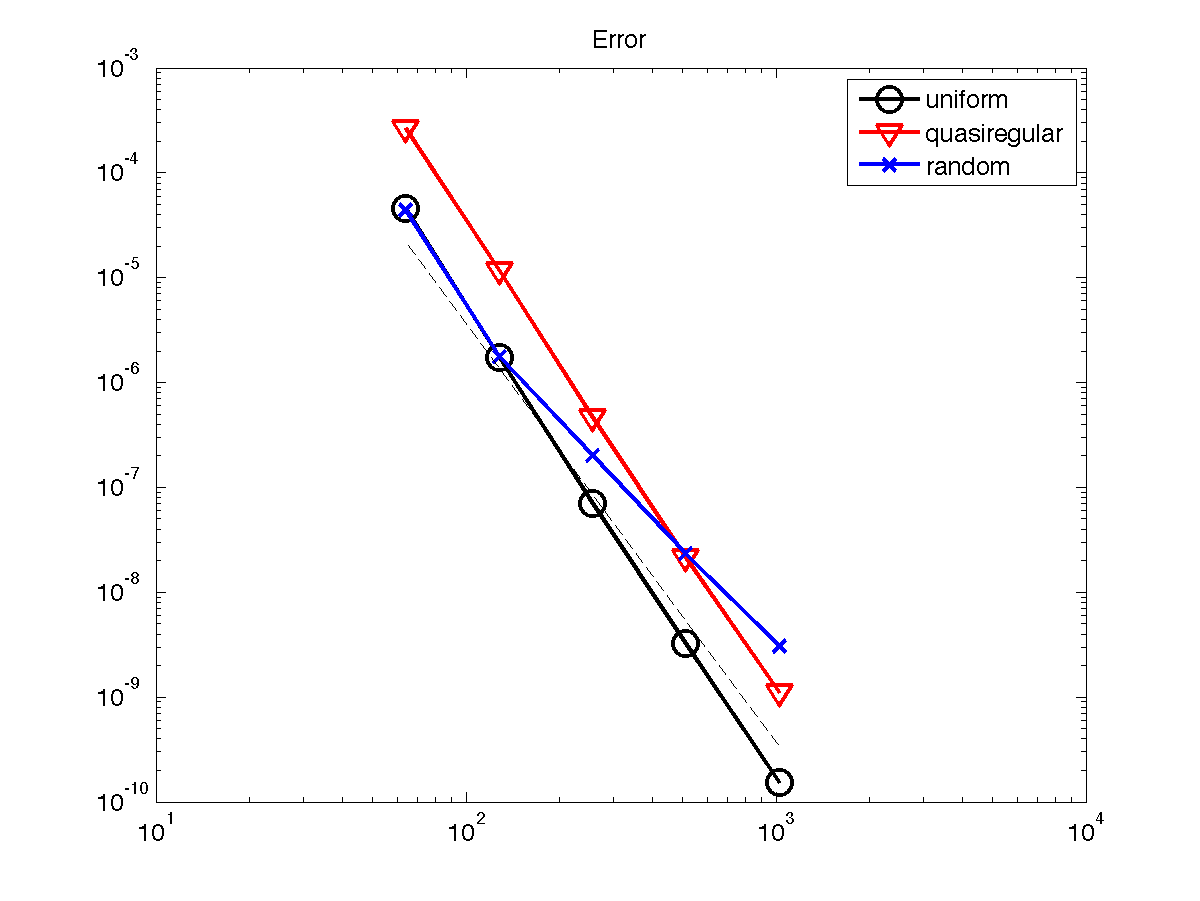}
\end{tabular}
\caption{Error decay under grid refinement for first (top-left), second (top-right), third (bottom-left) and fourth (bottom-right) order schemes. The dashed line indicates the expected decay in each case.}
\label{fig:convrate}
\end{figure}

\paragraph{Convergence tests}
Following \cite{XingShu:2005:WBSWEfd}, we compute the flow with initial data given by
\begin{equation}
\label{eq:test:Shu}
z(x)=\sin^2(\pi x)
\qquad
h(0,x) = 5+e^{\cos(2\pi x)}
\quad
q(0,x) = \sin(\cos(2\pi x))
\end{equation}
with periodic boundary conditions on the domain $[0,1]$. At time $t=0.1$ the solution is still smooth and we compare the numerical results with a reference solution computed with the fourth order scheme and $16384$ cells. The 1-norm of the errors appears in Figure \ref{fig:convrate} and the maximum entropy production is shown in Figure \ref{fig:entrate}  for all schemes and the three grid types considered.

All schemes have the expected accuracy, except for the fourth order scheme on the random grids, where the accuracy is slightly decreased due to the extreme irregularity of the grid. We point out however that, despite  the reduced decay rate, the actual values of the error of the fourth order scheme even on the random grid are orders of magnitude smaller than those obtained with the third order scheme with the same number of degrees of freedom.

\paragraph{Well-balancing}
We show a well-balancing test on the lake at rest solution using a bottom topography described by a uniformly distributed random variable sampled between $0$ and $1$, with water heigth at $h(x)+z(x)=1.5$. Table \ref{tab:wbtest} shows the well-balancing errors in the total water height and momentum, in the case of smooth nonuniform grids and random grids. Here $\Delta(h+z)_{j+1/2}= (h+z)_{j+1}- (h+z)_j$. All data are close to machine precision, as expected.

\begin{table}
\begin{center}
\begin{tabular}{|c|rrrr|rrrr|}
\hline
& \multicolumn{4}{c|}{$\|\Delta (h+z)\|_{\infty}$} & \multicolumn{4}{c|}{$\|q\|_{\infty}$} 
\\
\hline
Smooth & 100& 200 & 400 & 800 & 100& 200 & 400 & 800 \\
\hline
$p=1$ & 0 & 0 & 0 & 0 & 4.51e-16 &5.55e-16 & 5.00e-16 & 7.68e-16\\
$p=2$ & 0 & 2.22e-16 & 2.22e-16 & 2.22e-16 & 3.82e-16 &8.47e-16 & 7.36e-16 & 1.54e-15\\
$p=3$ & 0 &4.44e-16 & 4.44e-16 & 6.66e-16 & 6.87e-16 &1.47e-15 & 1.67e-15 & 2.47e-15\\
$p=4$ & 8.88e-16 &6.66e-16 & 1.55e-15 & 1.55e-15 & 9.89e-16 &1.82e-15 & 1.67e-15 & 1.90e-15\\
\hline
Random&\multicolumn{8}{c|}{}
\\
\hline
$p=1$ & 2.22e-16 & 2.22e-16 & 2.22e-16 & 2.22e-16 & 2.08e-16 &6.24e-16 & 6.77e-16 & 9.65e-16\\
$p=2$ & 2.22e-16 & 2.22e-16 & 2.22e-16 & 2.22e-16 & 2.91e-16 &7.25e-16 & 8.95e-16 & 9.99e-16\\
$p=3$ & 2.22e-16 & 6.66e-16 & 6.66e-16 & 6.66e-16 & 5.63e-16 &8.47e-16 & 9.94e-16 & 1.28e-15\\
$p=4$ & 6.66e-16 & 8.88e-16 & 1.33e-15 & 1.11e-15 & 8.68e-16 &7.94e-16 & 1.11e-15 & 1.43e-15\\
\hline
\end{tabular}
\end{center}
\caption{Lake at rest test: well-balancing errors with rough bottom. }
\label{tab:wbtest}
\end{table}

\paragraph{Small perturbation of a lake at rest}
The domain is $x\in[0,2]$, the bottom and initial total height are given by
\begin{equation}
\label{eq:smallpulse}
z(x)=\begin{cases}
0.25(1+\cos(10\pi(x-0.5))) & 1.2\leq x\leq 1.4\\
0 &\text{otherwise}
\end{cases}
\qquad
H(x,0)=1+ 0.001 \chi_{[1.1,1.2]}(x)
\end{equation}

\begin{figure}
\includegraphics[width=0.45\linewidth]{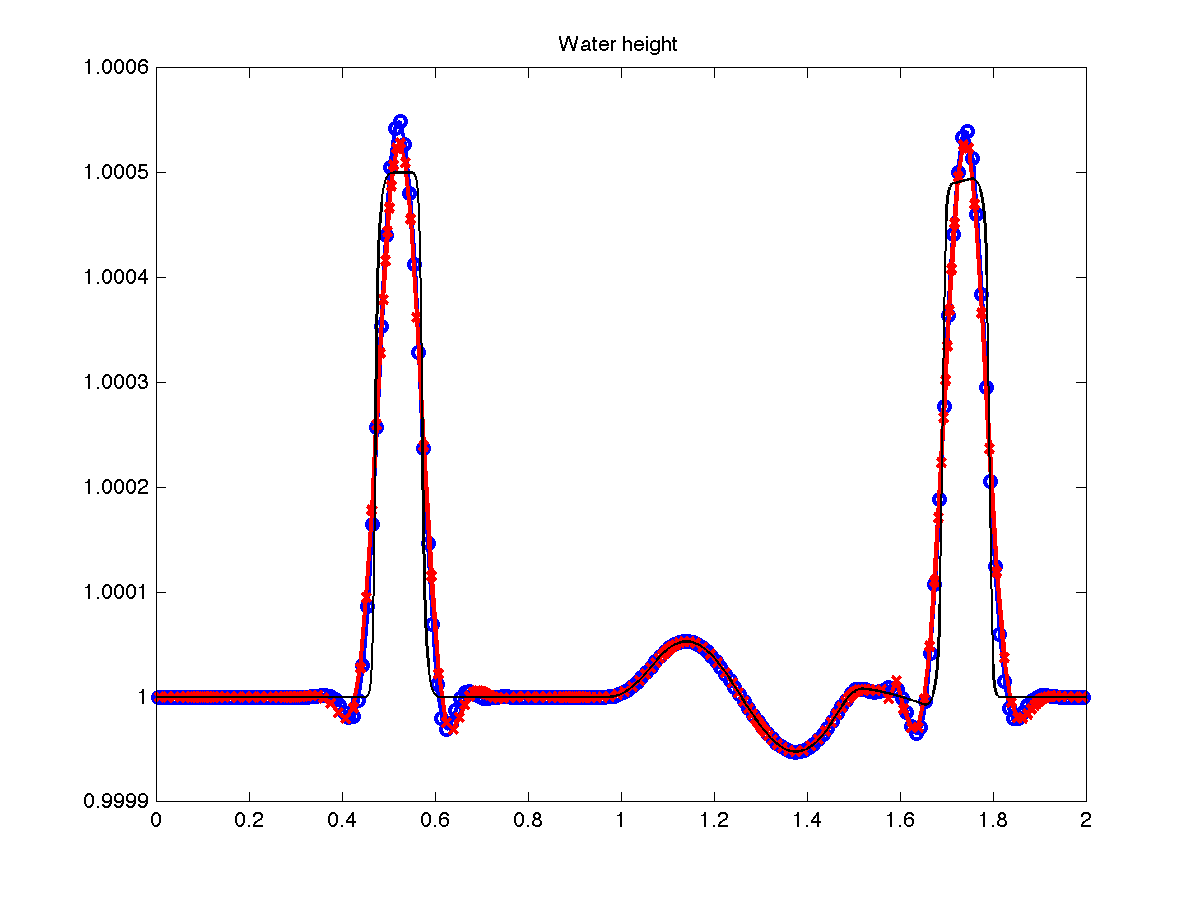}
\hfill
\includegraphics[width=0.45\linewidth]{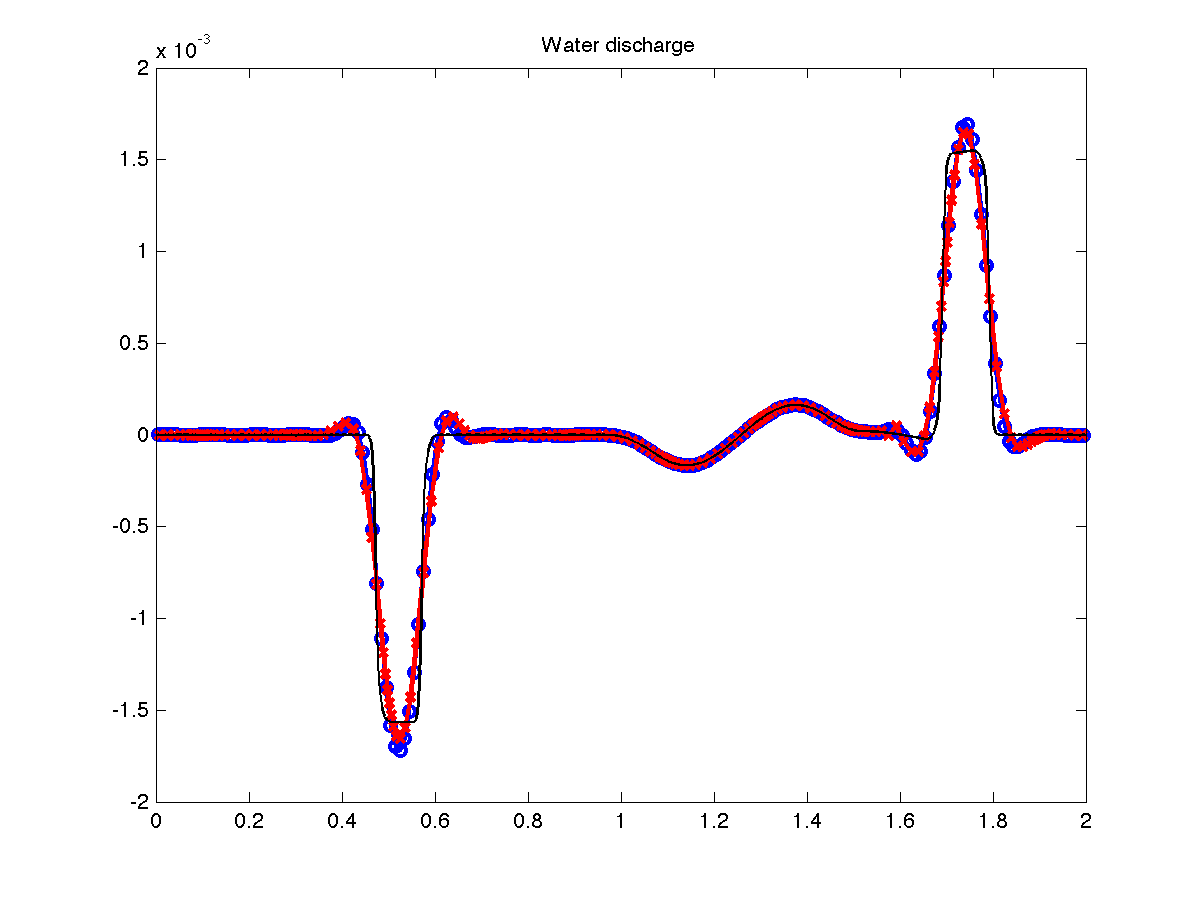}
\caption{LeVeque's test \eqref{eq:smallpulse}. Third order scheme on a uniform (blue circles) and {\em quasi-regular} grid (red crosses)  on top of a reference solution (black solid line).}
\label{fig:smallpulse:3:regular}
\end{figure}

\begin{figure}
\includegraphics[width=0.45\linewidth]{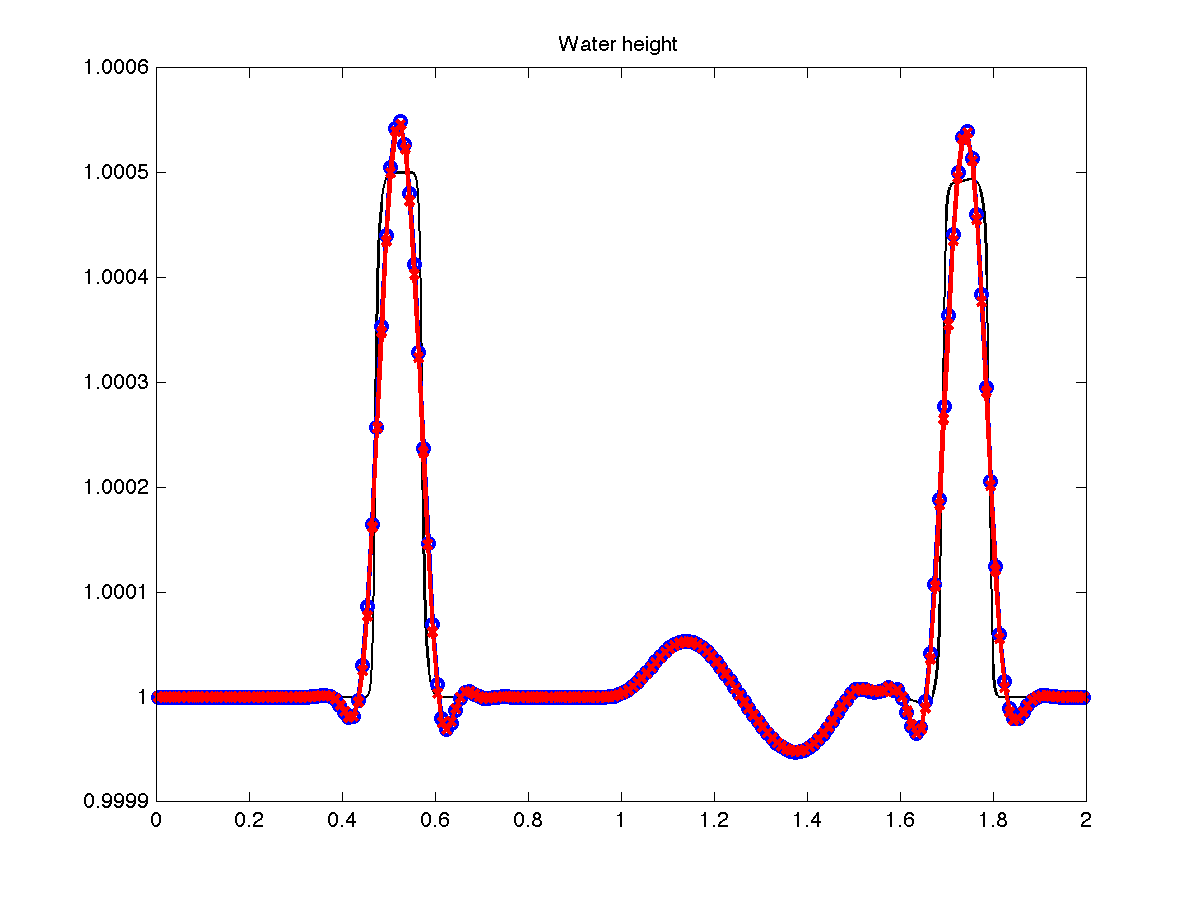}
\hfill
\includegraphics[width=0.45\linewidth]{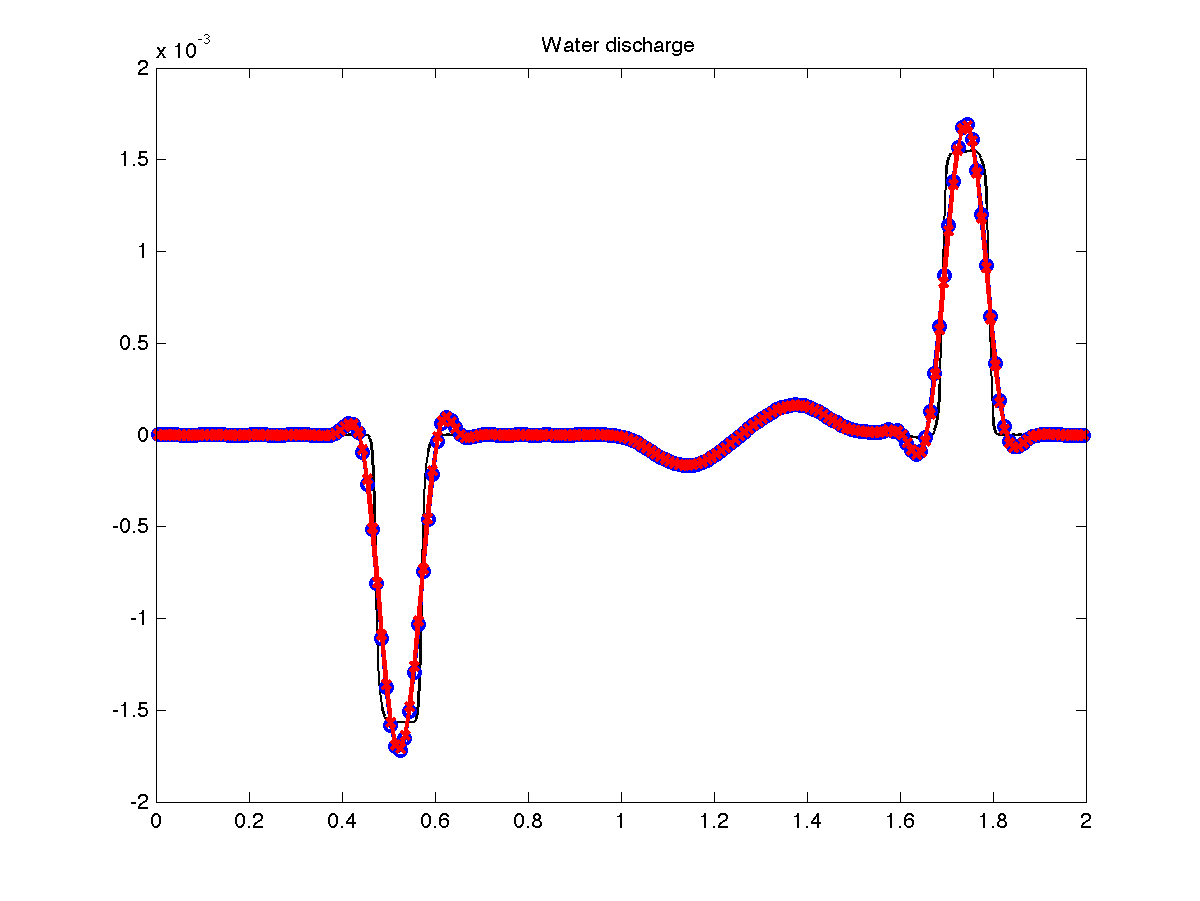}
\caption{LeVeque's test \eqref{eq:smallpulse}. Third order scheme on uniform (blue circles) and {\em random} grid (red crosses) on top of a reference solution (black solid line)}
\label{fig:smallpulse:3:random}
\end{figure}

\begin{figure}
\includegraphics[width=0.45\linewidth]{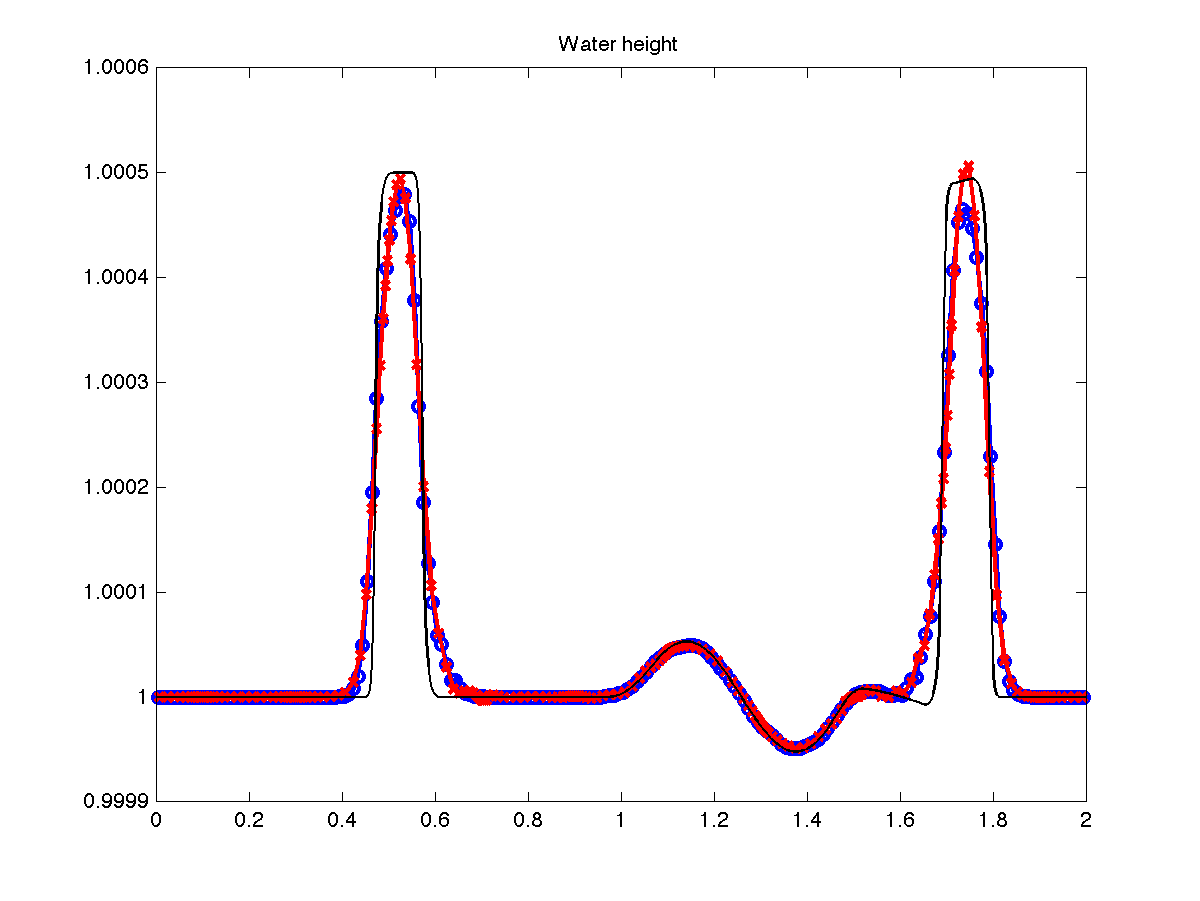}
\hfill
\includegraphics[width=0.45\linewidth]{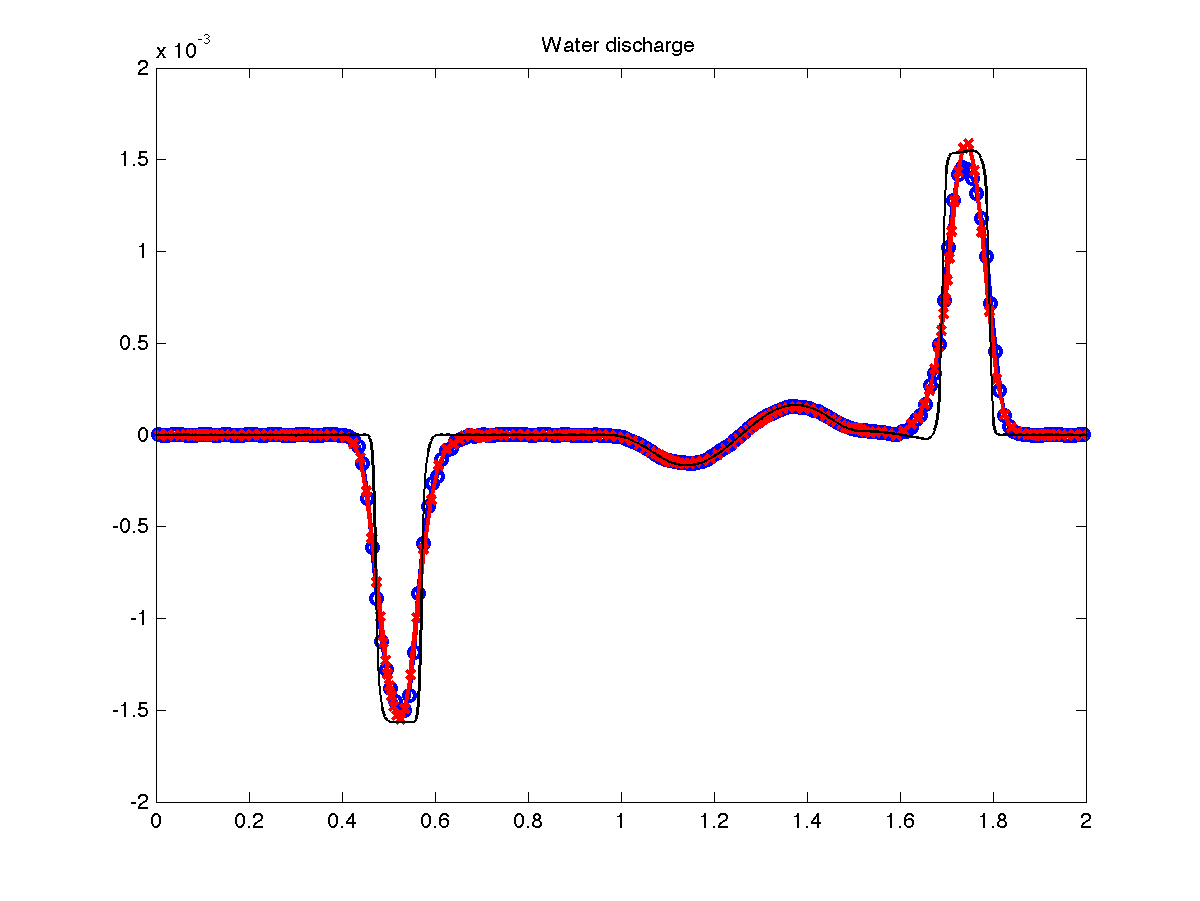}
\caption{LeVeque's test \eqref{eq:smallpulse}. Fourth order scheme on  uniform (blue circles) and {\em quasi-regular} grid (red crosses), on top of a reference solution (black solid line)}
\label{fig:smallpulse:4:regular}
\end{figure}

\begin{figure}
\includegraphics[width=0.45\linewidth]{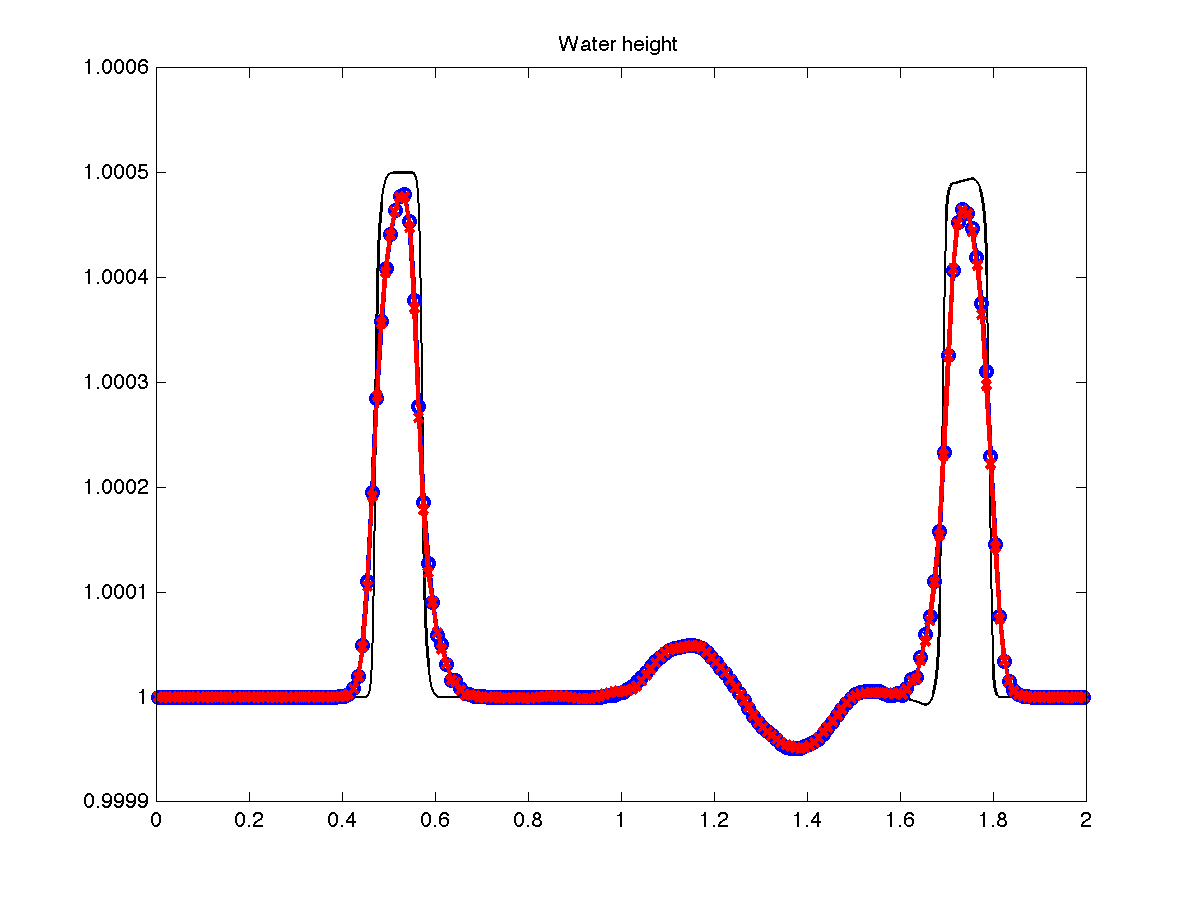}
\hfill
\includegraphics[width=0.45\linewidth]{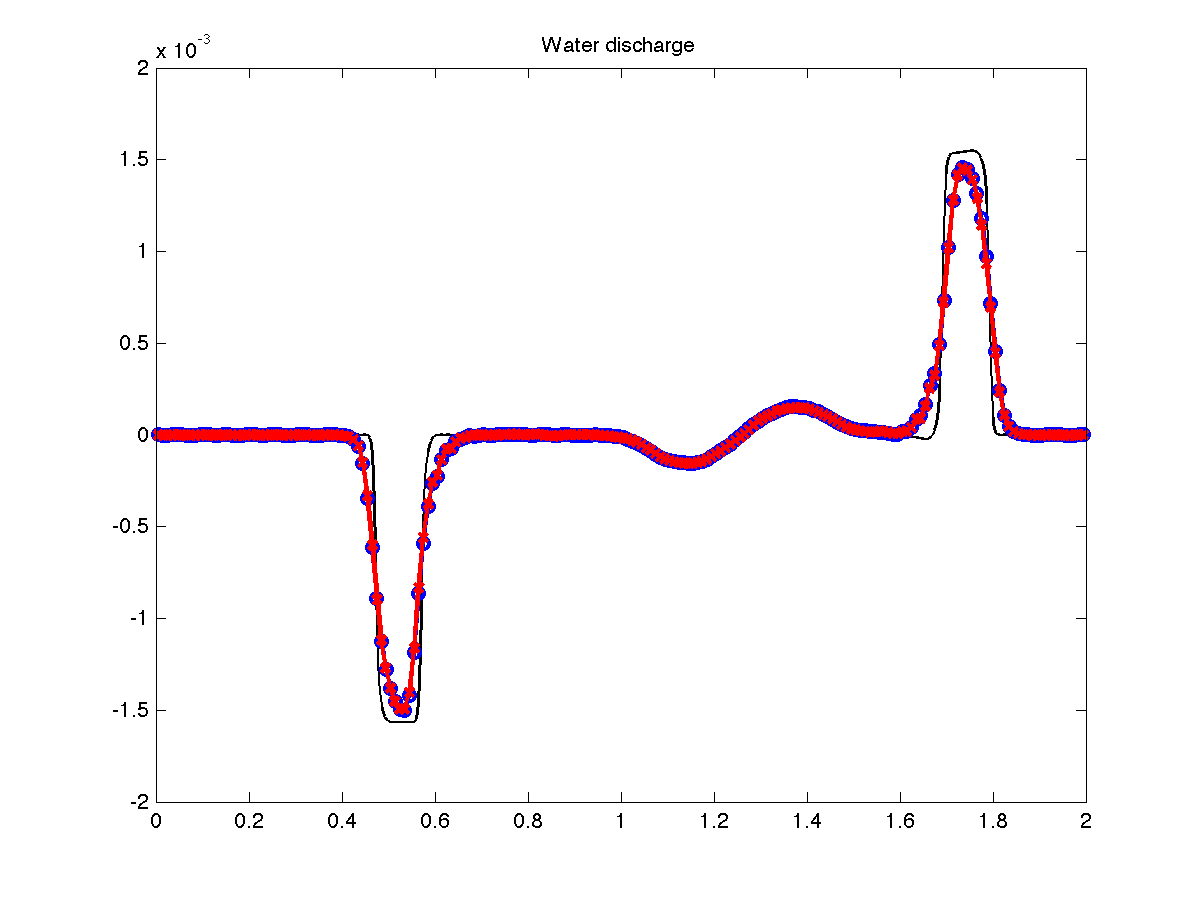}
\caption{LeVeque's test \eqref{eq:smallpulse}. Fourth order scheme on uniform (blue circles) and {\em random} grids (red crosses), on top of a reference solution (black solid line).}
\label{fig:smallpulse:4:random}
\end{figure}

This test was first used by LeVeque in \cite{LeVeque} with a second order scheme, but here we use it with a smaller perturbation for the third and fourth order schemes, as in \cite{NatvigEtAl}. This test requires a well-balanced scheme to resolve correctly the small perturbations which otherwise would be hidden by numerical noise. The solutions are shown in Fig \ref{fig:smallpulse:3:regular} and \ref{fig:smallpulse:3:random} for the third order scheme and Fig \ref{fig:smallpulse:4:regular} and \ref{fig:smallpulse:4:random} for the fourth order one. In each of the figures the numerical solution obtained with the uniform grid is compared with the one obtained on a non-uniform mesh. It can be seen that the pulse is well-resolved in all cases and the results obtained with a uniform grid can be perfectly superposed on those computed with the uniform ones. In this test, the parameter $\epsilon$ in the nonlinear weights of the WENO schemes is set to $10^{-12}$, as pointed out in \cite{NatvigEtAl}.

\paragraph{Moving water equilibria}
Since our schemes are well-balanced around the lake-at-rest equilibrium, one does not expect them to compute moving water equilibria at machine precision. Here we show two tests. In the first case we consider a transcritical steady state with a shock, over the parabolic hump
\[ 
z(x)=\begin{cases}
(0.2-0.05*(x-10)^2) & 8\leq x \leq 12\\
0 &\text{otherwise}
\end{cases}
\]
in the domain $[0,25]$. We consider the steady state solution with $q(x)=0.18$, with Dirichlet boundary conditions $q=0.18$ at $x=0$ and $h=0.33$ at $x=25$. The solution has a steady shock at $x=11.665504281554291$. The computation was initialized with the exact steady state solution (see for example the Appendix A of \cite{Karni}) and the numerical integration was performed until $t=50$.

\begin{figure}
\includegraphics[width=0.45\linewidth]{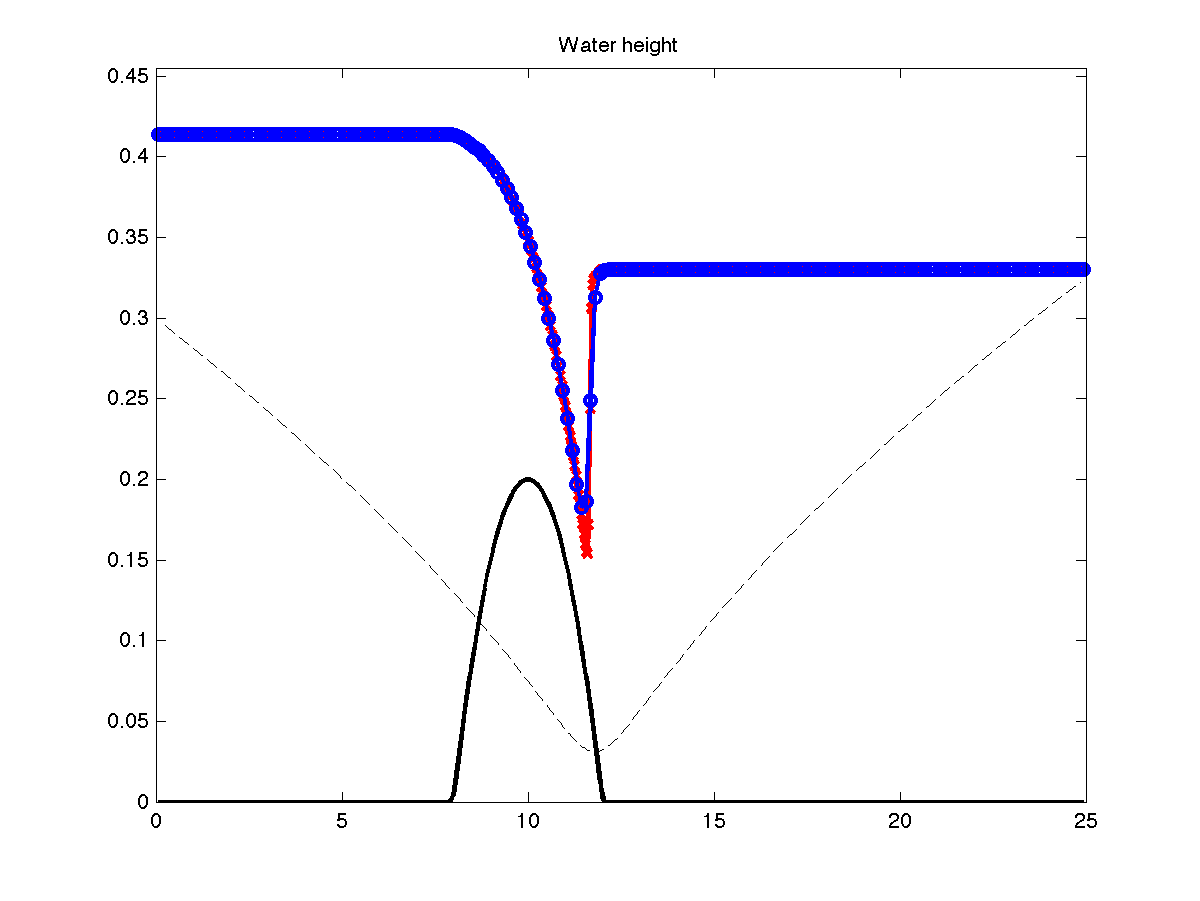}
\hfill
\includegraphics[width=0.45\linewidth]{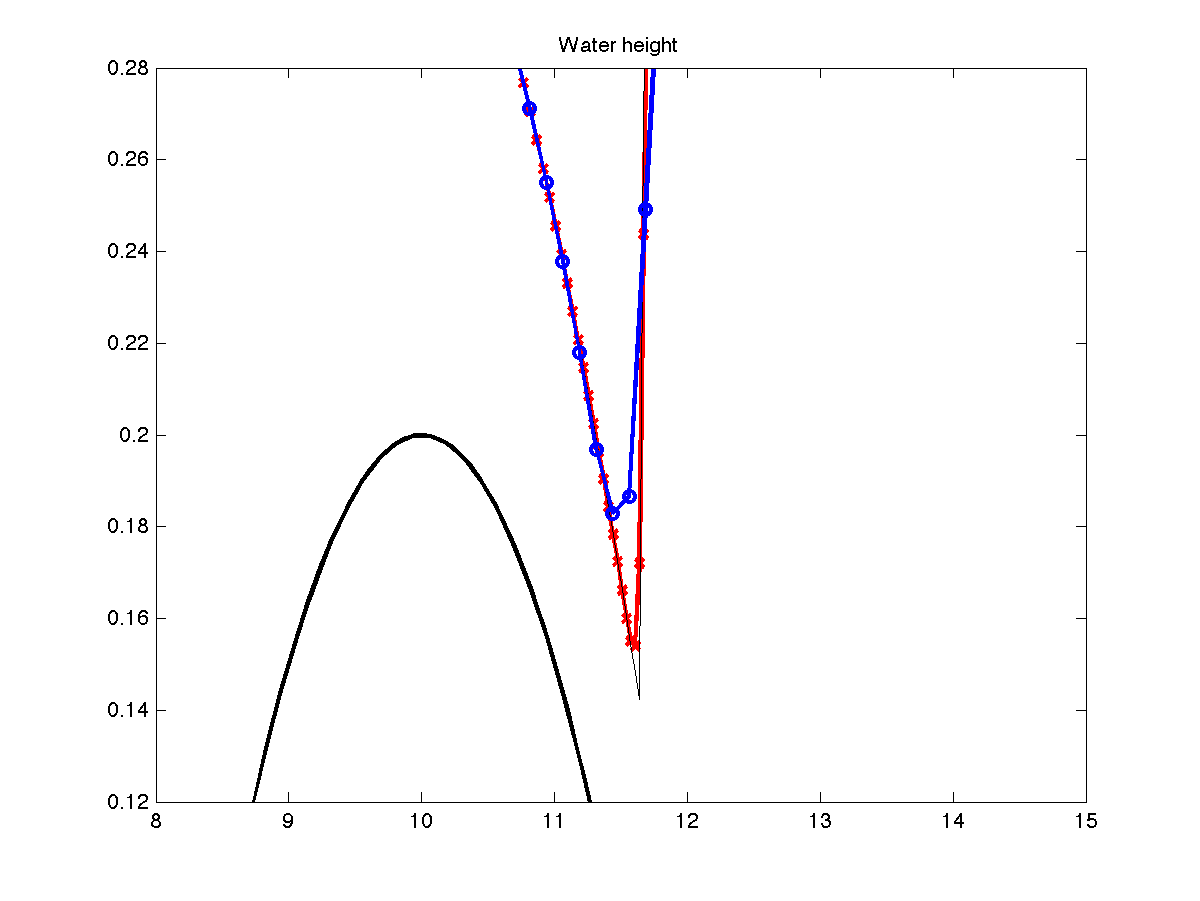}
\caption{Steady solution with transcritical shock, approximated with a third order scheme (uniform and adapted grids). The dashed line in the left panel is the local grid size in the non-uniform grid.}
\label{fig:steady:transshock:3}
\end{figure}

\begin{figure}
\includegraphics[width=0.45\linewidth]{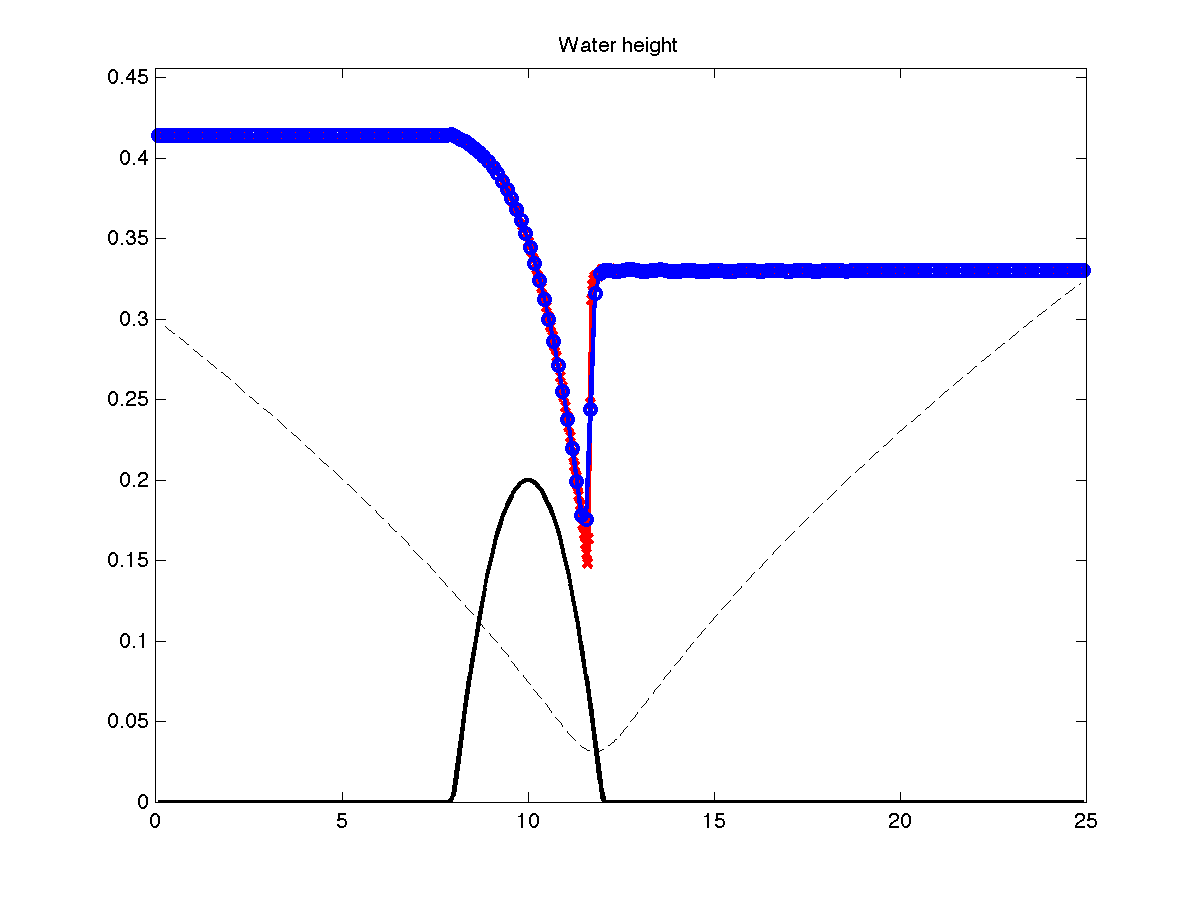}
\hfill
\includegraphics[width=0.45\linewidth]{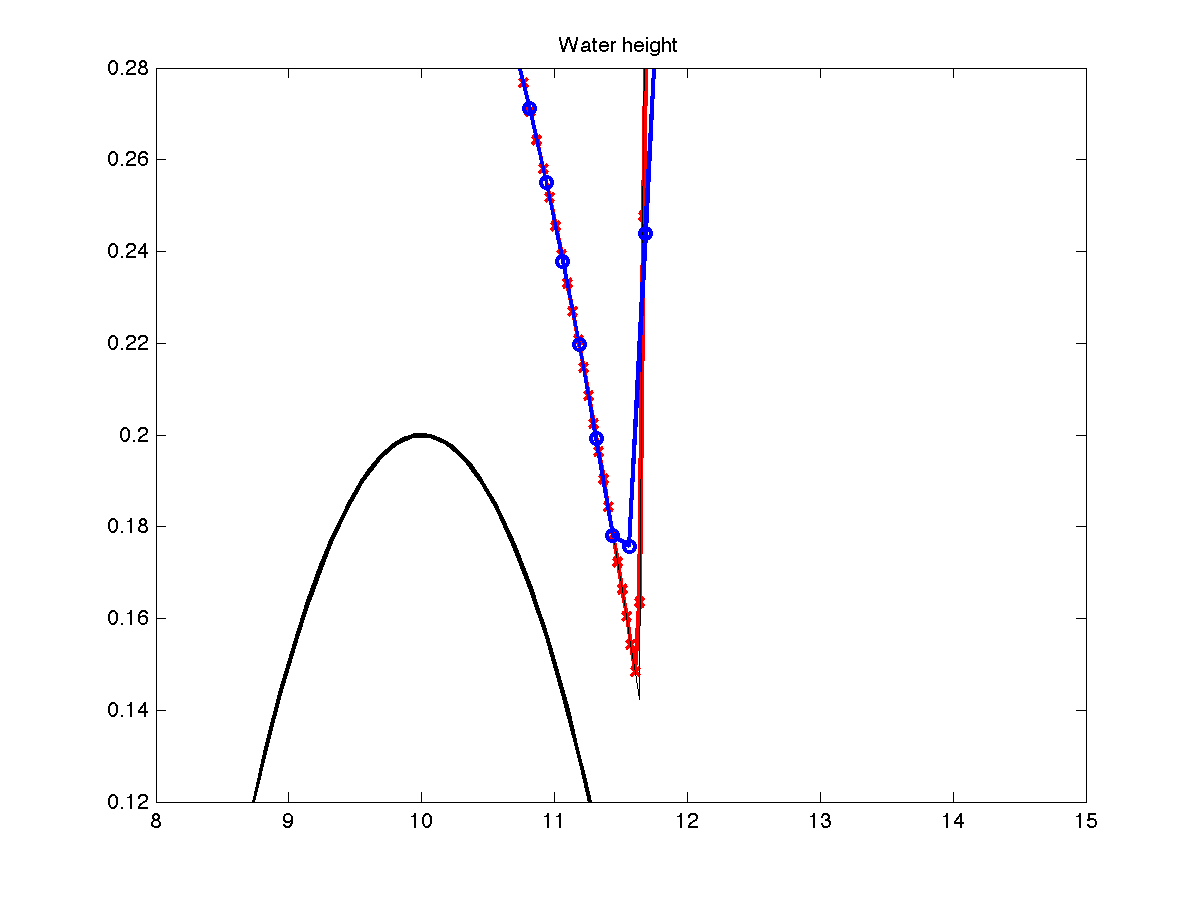}
\caption{Steady solution with transcritical shock, approximated with a fourth order scheme (uniform and adapted grids). The dashed line in the left panel is the local grid size in the non-uniform grid.}
\label{fig:steady:transshock:4}
\end{figure}


\begin{table}
\begin{center}
\begin{tabular}{|c|rr|rr|rr|rr|}
\hline
& \multicolumn{2}{c|}{$p=1$} 
& \multicolumn{2}{c|}{$p=2$} 
& \multicolumn{2}{c|}{$p=3$} 
& \multicolumn{2}{c|}{$p=4$} 
\\
\hline
Uniform & error& rate &  error& rate & error& rate & error& rate  \\
\hline
$100$ & 1.96e-1 & --   & 5.54e-2 &   -- & 2.02e-2 &   -- & 2.92e-3 &   --\\
$200$ & 1.17e-1 & 0.74 & 1.42e-2 & 1.96 & 4.26e-3 & 2.24 & 1.40e-4 & 4.38\\
$400$ & 6.35e-2 & 0.89 & 3.29e-3 & 2.11 & 4.87e-4 & 3.13 & 5.12e-6 & 4.77\\
$800$ & 3.26e-2 & 0.96 & 8.08e-4 & 2.03 & 3.89e-5 & 3.65 & 1.60e-7 & 5.00\\
\hline
Adapted&\multicolumn{8}{c|}{} \\ 
\hline
$100$ & 9.20e-2 &   -- & 6.96e-3 & --   & 9.78e-4 &   -- & 4.54e-5 & --      \\
$200$ & 4.67e-2 & 0.97 & 1.71e-3 & 2.02 & 7.97e-5 & 3.62 & 1.36e-6 & 5.07 \\
$400$ & 2.34e-2 & 0.99 & 4.25e-4 & 2.01 & 6.57e-6 & 3.60 & 3.87e-8 & 5.13\\
$800$ & 1.17e-2 & 1.00 & 1.06e-4 & 2.01 & 5.63e-7 & 3.55 & 1.25e-9 & 4.95 \\
\hline
\end{tabular}
\end{center}
\caption{Well-balancing errors for the subcritical steady state with gaussian bottom.}
\label{tab:wbsubcritial}
\end{table}

We show the solutions computed with uniform grids and with a grid refined ad-hoc around the shock position (see Eq \eqref{eq:locallyrefined}) with the scheme of order three (Figures \ref{fig:steady:transshock:3}) and four (Figure \ref{fig:steady:transshock:4}). The figures report with a dashed line the local cell size of the nonuniform grid, which is refined close to the shock. The right panels of each figure show a zoom on the shock and it is clear that the adapted solution (in red with crosses) approximates better the exact solution (thin black line) than the solution obtained with a uniform grid with the same numer of points (blue line with dots), with no spurious oscillations.


In order to quantify the improvement due to the adapted grid and the rate of convergence of the schemes on moving water equilibria, we consider a smooth test problem, namely a subcritical steady flow over the smooth bump $z(x)=0.2e^{-(x-12.5)^2}$ on the domain $[0,25]$. The numerical scheme was initialized with the exact solution and the flow computed until $t=10$. Since the behaviour of the errors on the water height and on  momentum is very similar, only the former are reported in Table \ref{tab:wbsubcritial}. The first and second order schemes show the expected rates of convergence, while the third and fourth order ones have convergence rates well above the expected values (respectively $3.60$ and $5.00$). 

We also consider nonuniform grids that are finer on the hump and coarser on the flat portion of the bottom function, namely those given by Eq. \eqref{eq:locallyrefined} with $\overline{w}=12.5/25=0.5$. The errors on the adapted grids are much smaller than the corresponding uniform grids and the convergence rates are confirmed also on nonuniform grids.


\subsection{Numerical entropy production}

\begin{figure}
\begin{tabular}{cc}
First order
&
Second order
\\
\includegraphics[width=0.45\linewidth]{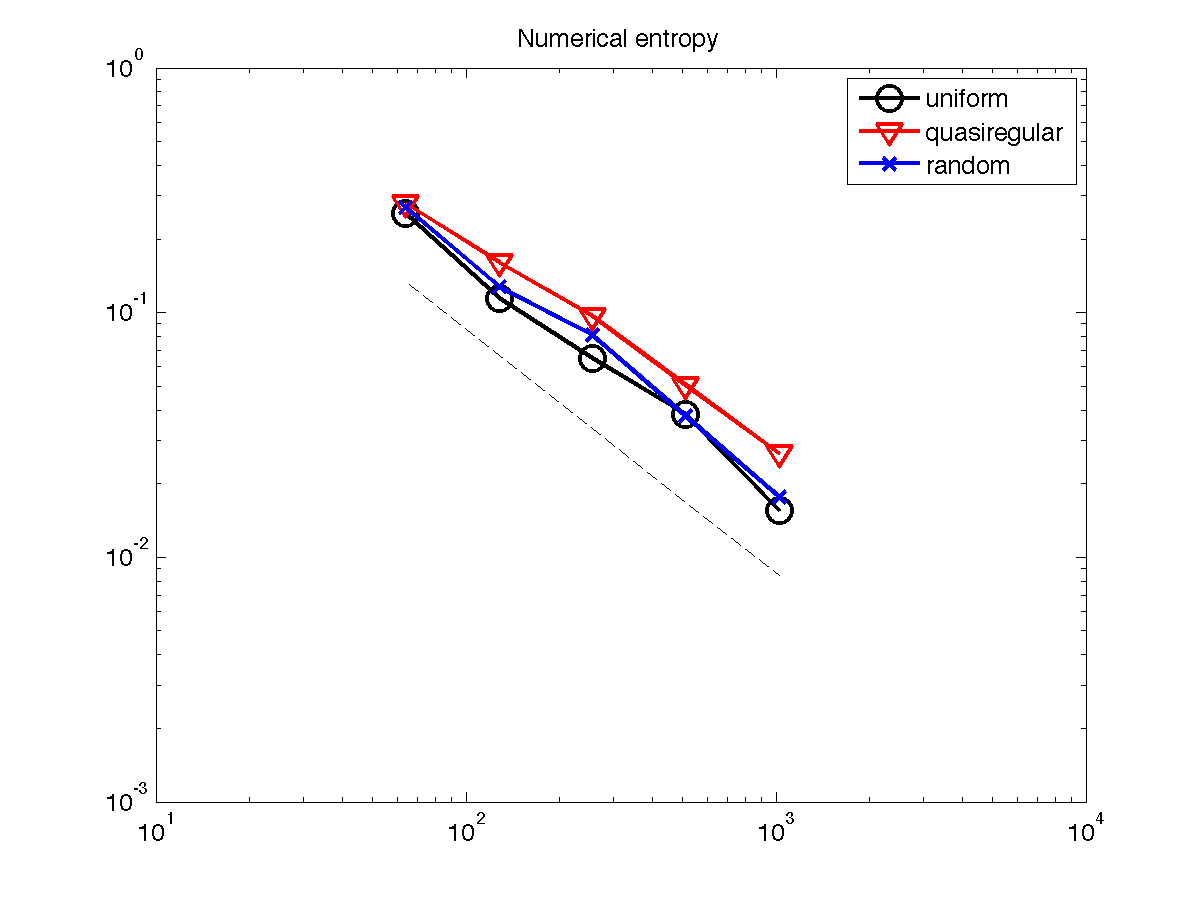}
&
\includegraphics[width=0.45\linewidth]{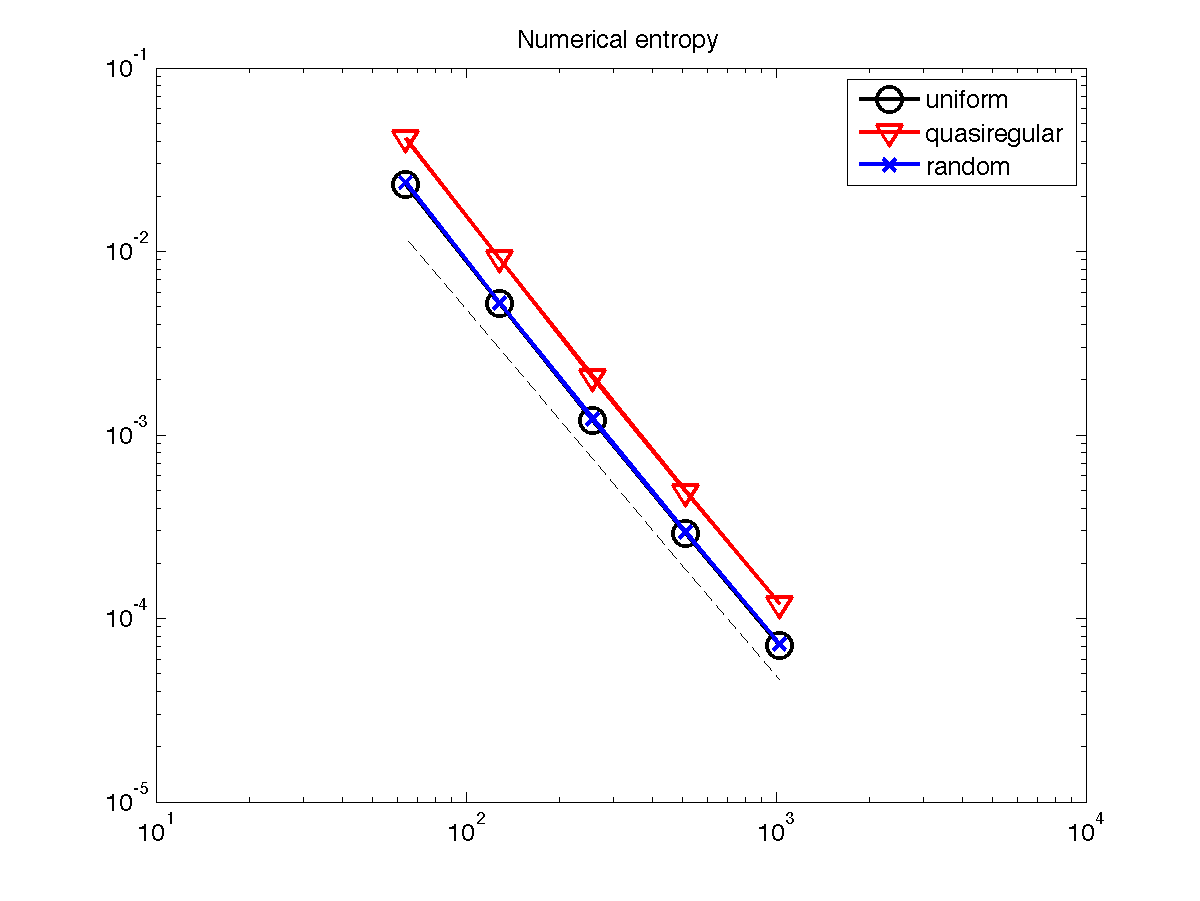}
\\
Third order
&
Fourth order
\\
\includegraphics[width=0.45\linewidth]{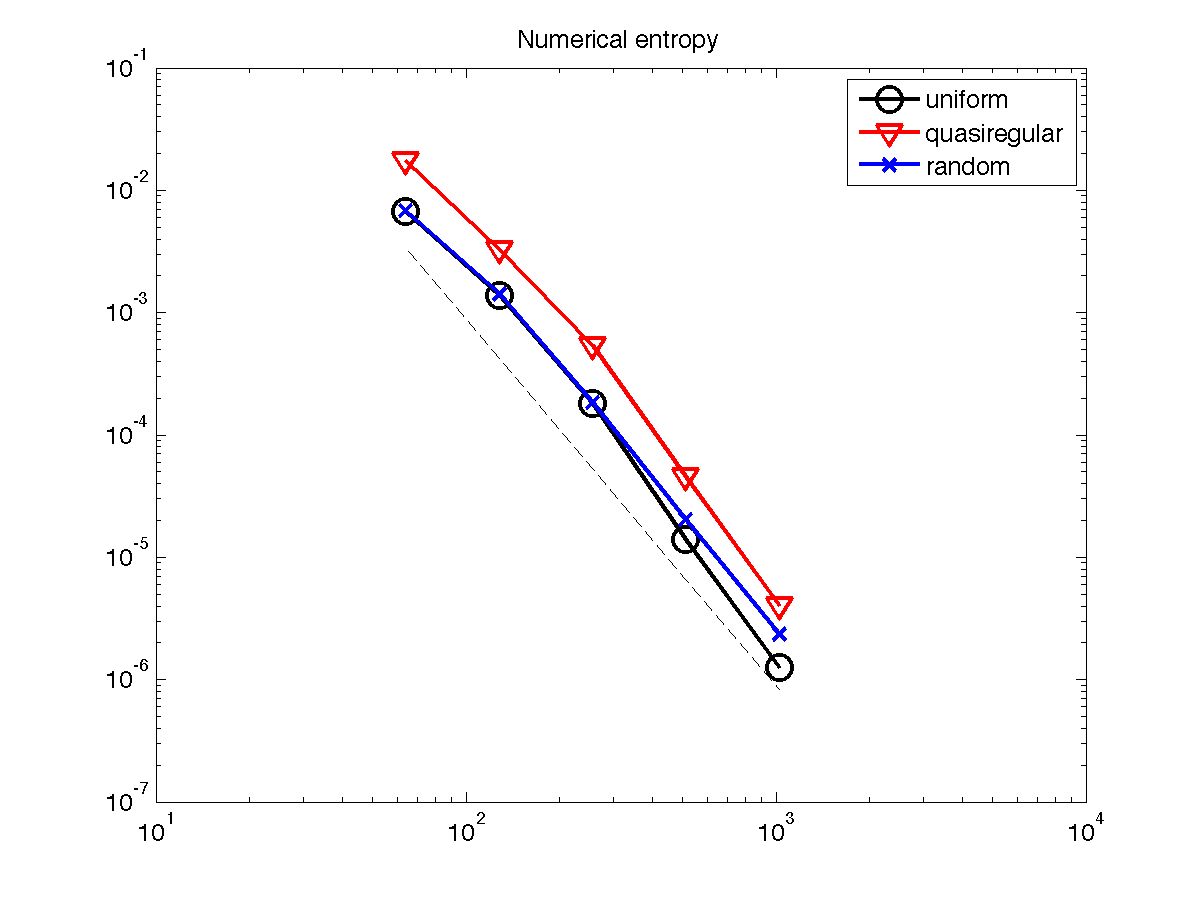}
&
\includegraphics[width=0.45\linewidth]{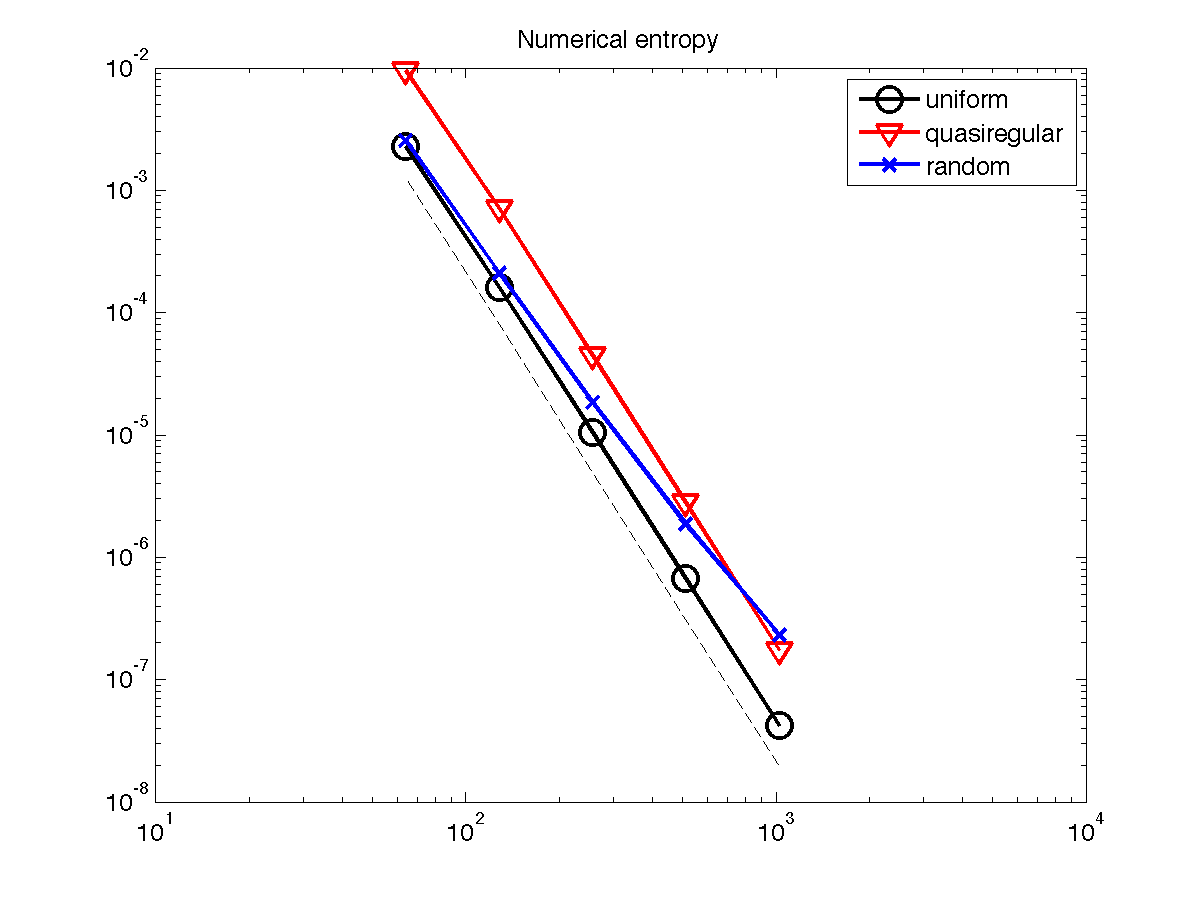}
\end{tabular}
\caption{Numerical entropy production decay under grid refinement for first (top-left), second (top-right), third (bottom-left) and fourth (bottom-right) order schemes. The dashed line indicates the expected decay in each case.}
\label{fig:entrate}
\end{figure}

\paragraph{Rate of decay on smooth flows.} 
Figure \ref{fig:entrate} shows the numerical entropy production in the smooth test \eqref{eq:test:Shu} on several grid types. It is apparent that the decay rate, as expected, follows the order of accuracy of the corresponding schemes. Moreover, comparing this figure with Figure \ref{fig:convrate}, we note that the entropy decay mimics exactly the behaviour of the error, even in the case of the slight deterioration of accuracy observed on the random grid for the fourth order scheme.

\begin{figure}
\begin{center}
\begin{tabular}{cc}
\includegraphics[width=0.45\linewidth]{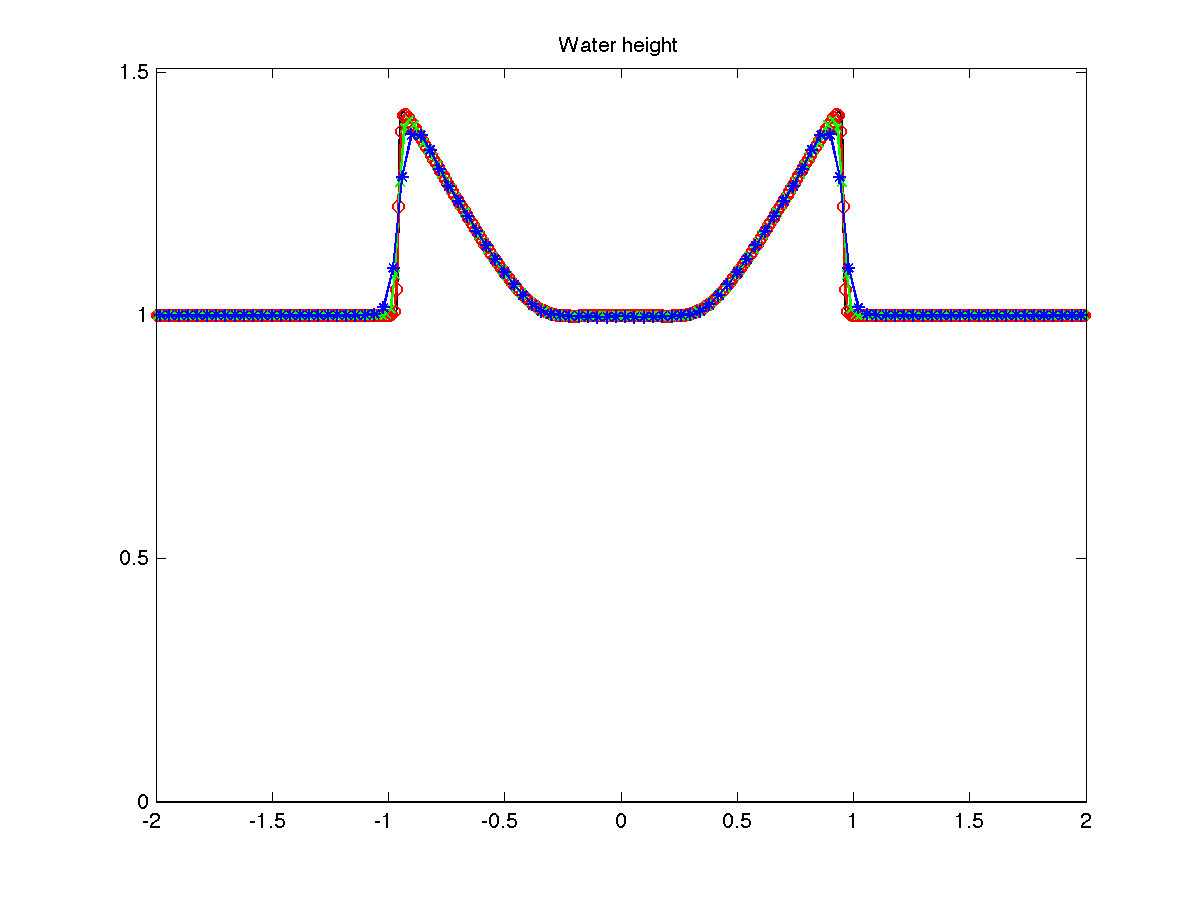}
&
\includegraphics[width=0.45\linewidth]{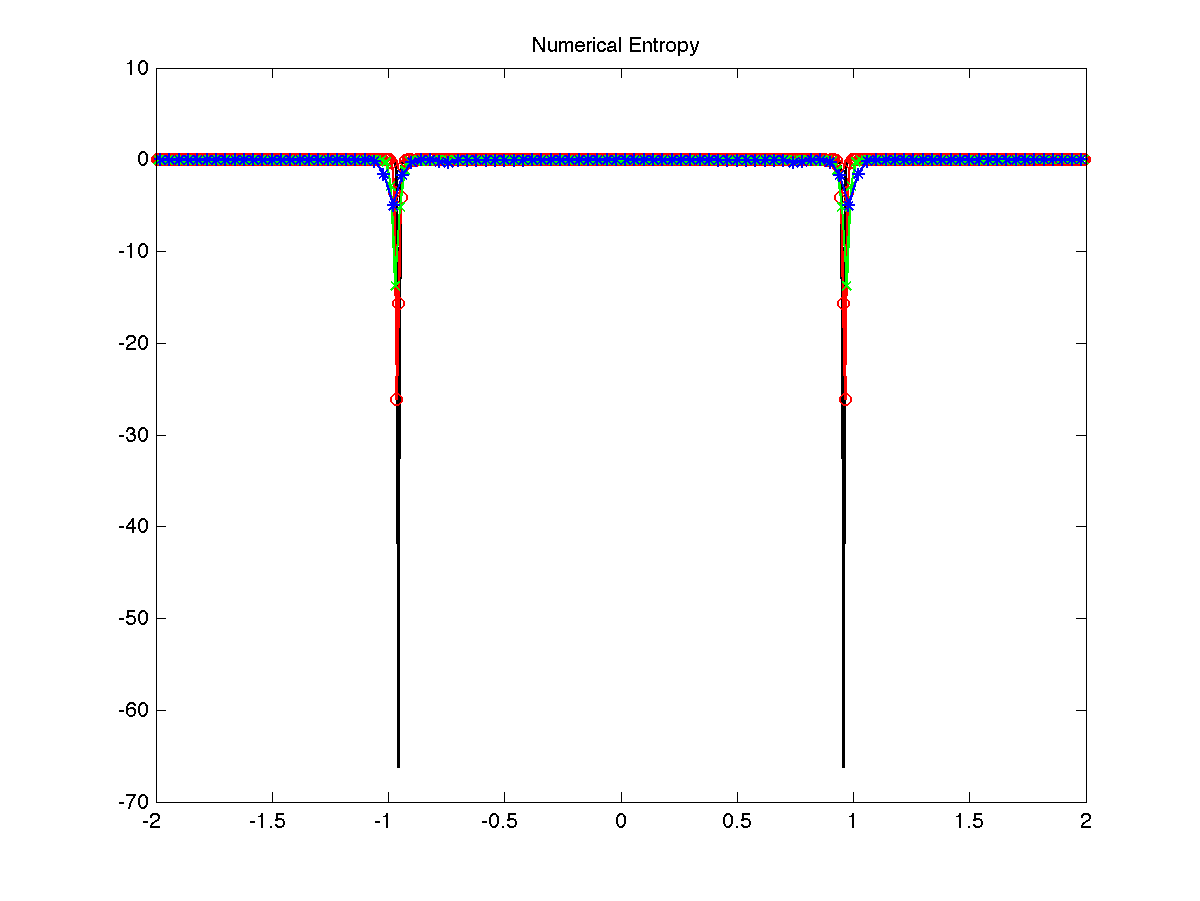}
\\
\includegraphics[width=0.45\linewidth]{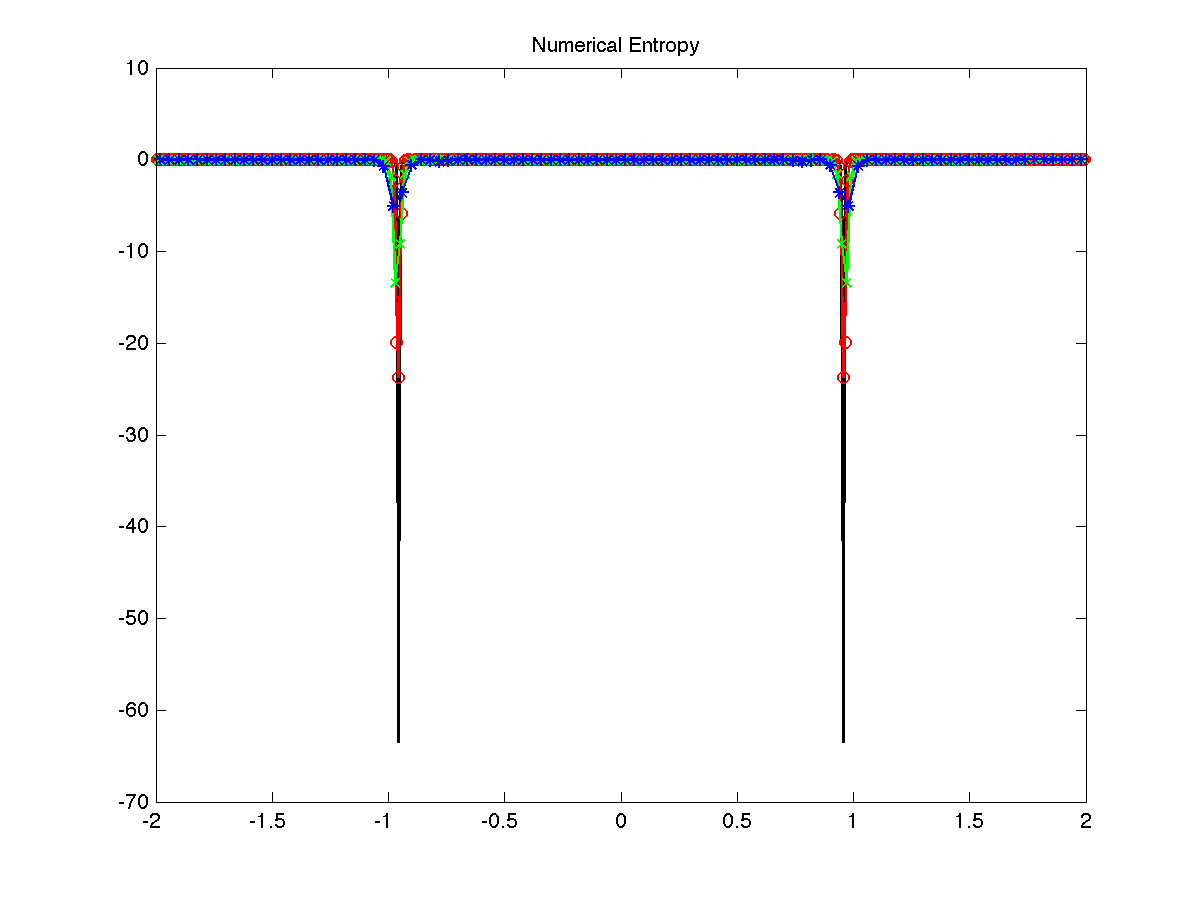}
&
\includegraphics[width=0.45\linewidth]{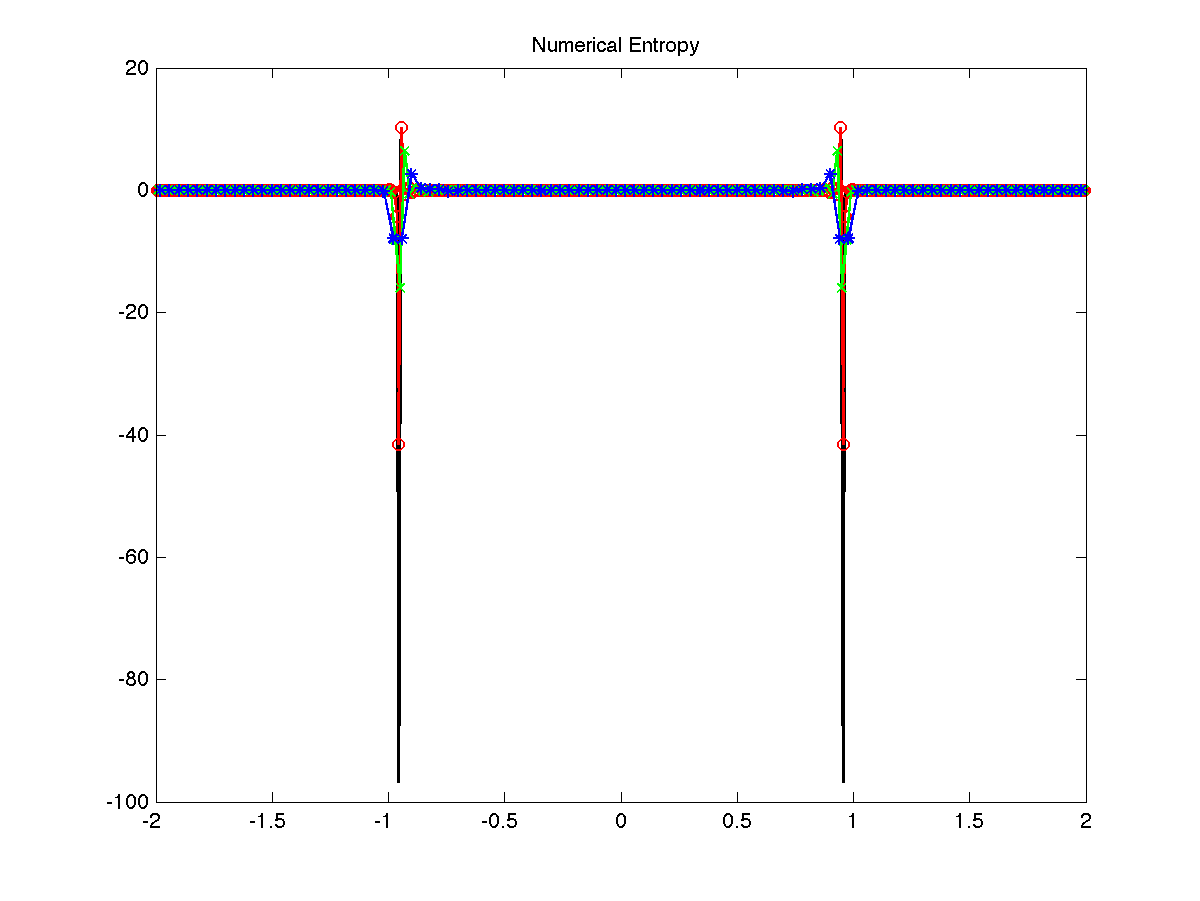}
\end{tabular}
\end{center}
\caption{Entropy production on shocks under grid refinement for several schemes. Top-left: water height. Top-right: second oder scheme. Bottom-left: third order scheme. Bottom right: fourth order scheme. $N=800$  (black solid line), $N=400$  (red line with circles), $N=200$  (green line with crosses), $N=100$ (blue line with stars). }
\label{fig:stonato}
\end{figure}

\paragraph{Two shocks.}
We set up initial data with a flat bottom, water at rest and 
$h(0,x) = e^{-50x^2}$ on the domain $[-2,2]$. As the flow evolves, two shocks form and separate from each other: at $t=0.2$ the computed water height is depicted in the top-left plot of Figure \ref{fig:stonato}. Each of the other panels of Figure \ref{fig:stonato} shows the entropy residual obtained with four different grid sizes. The results for second, third and fourth order schemes appear in the top-right, lower left and lower right panels respectively.  In all three cases it can be seen that the numerical entropy production on the two shocks increases under grid refinement like $1/h$. On the other hand, the magnitude of the peak of the numerical entropy production does not depend on the order of the scheme. This is to be contrasted with the numerical entropy production on smooth flows just shown, where one observes entropy residuals of $O(h^p)$, where $p$ is the order of the scheme.

Due to the different orders of magnitude of the numerical entropy production in the smooth regions of the flows and around shocks, it can be concluded that the entropy residual provides an effective discontinuity detector, expecially in the case of high order schemes.

\begin{figure}
\includegraphics[width=0.45\linewidth]{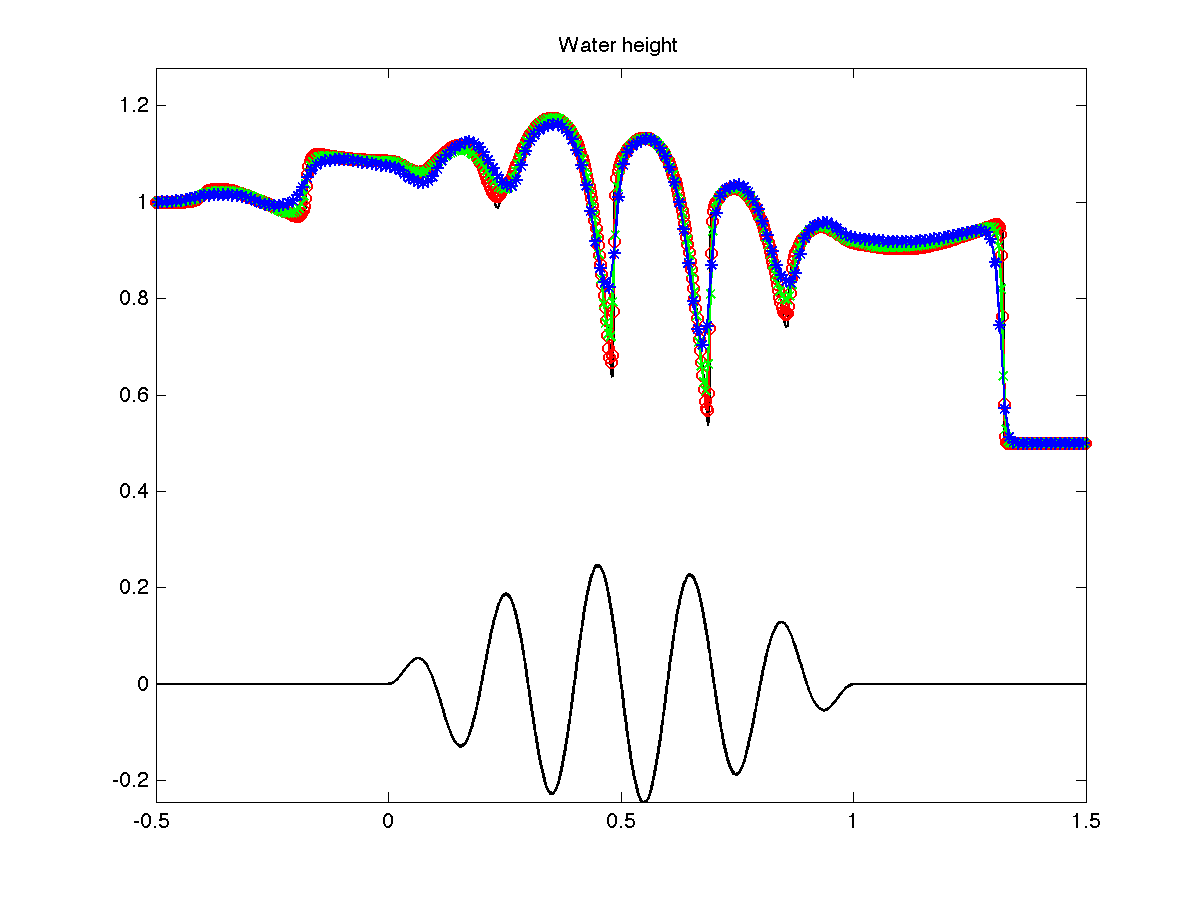}
\hfill
\includegraphics[width=0.45\linewidth]{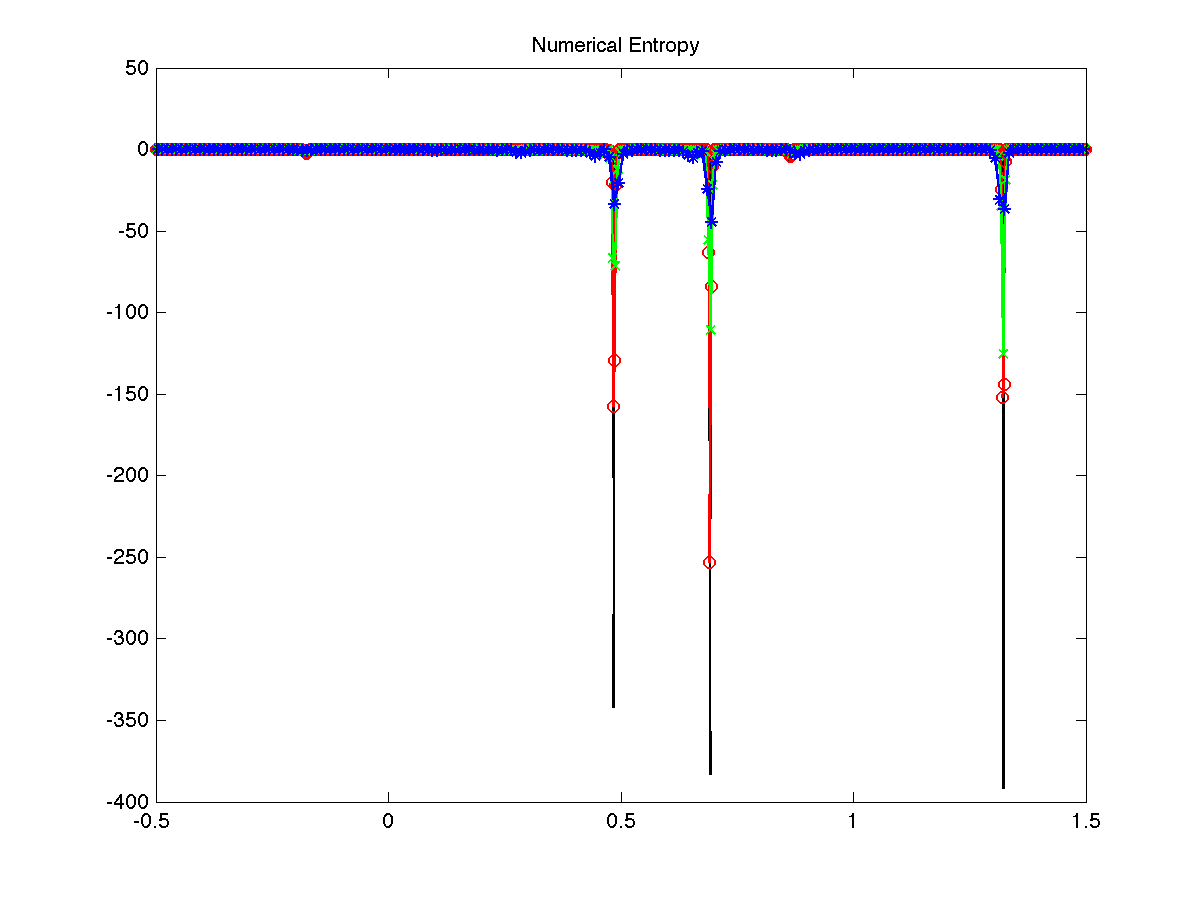}
\caption{Stream on artificial river bed. Left: water height. Right: numerical entropy production. $N=1600$ (black solid line), $N=800$ (red line with circles), $N=400$  (green line with crosses) and $N=200$ points (blue line with stars).}
\label{fig:sincanal}
\end{figure}

\paragraph{Stream on artificial river bed.}
In the domain $[-.5,1.5]$ we consider the bottom topography and initial conditions:
\begin{eqnarray}
\label{eq:SinCanal}
& z(x) = 
\begin{cases}
 \sin(10\pi x) x(1-x) & x\in[0,1]\\
 0& \text{otherwise}
\end{cases}
\\
& H(0,x)= 
\begin{cases}
 1.0 & x< -0.2\\
 0.5 & x\geq -0.2
\end{cases}
\quad
q(0,x)=
\begin{cases}
 \tfrac12\sqrt{\tfrac32 g} & x< -0.2\\
 0.0 & x\geq -0.2
\end{cases} \nonumber
\end{eqnarray}
We integrate with free flow boundary conditions until $t=0.4$, when the shock originated from the Riemann problem has overcome the irregularity in the bottom topography (see the left panel of Figure \ref{fig:sincanal}). The right panel compares the numerical entropy production of the second order scheme with grid size from $200$ to $1600$. The peaks in the numerical  entropy production clearly show the location of the shocks and have the expected $O(1/h)$ behaviour.

\begin{figure}
\begin{center}
\begin{tabular}{cc}
First order
&
Second order
\\
\includegraphics[width=0.45\linewidth]{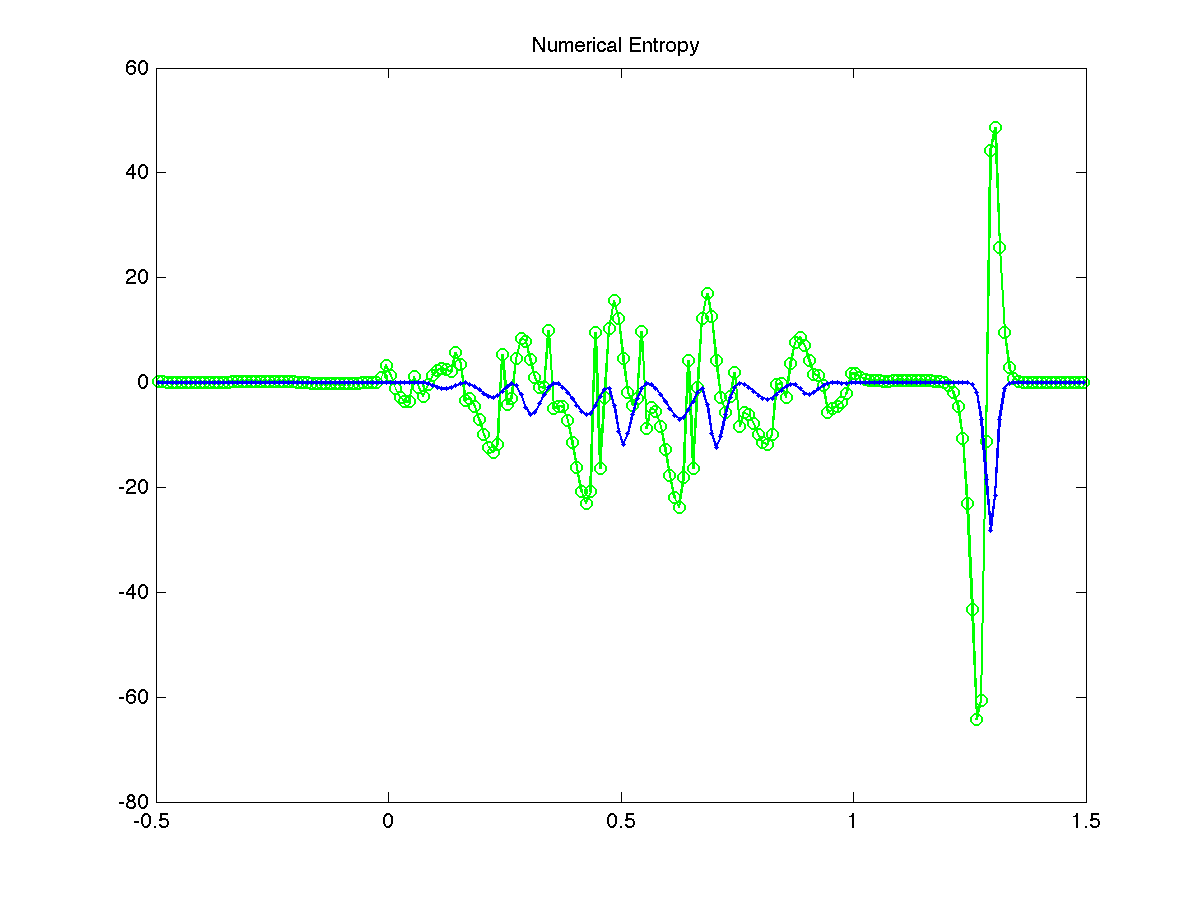}
&
\includegraphics[width=0.45\linewidth]{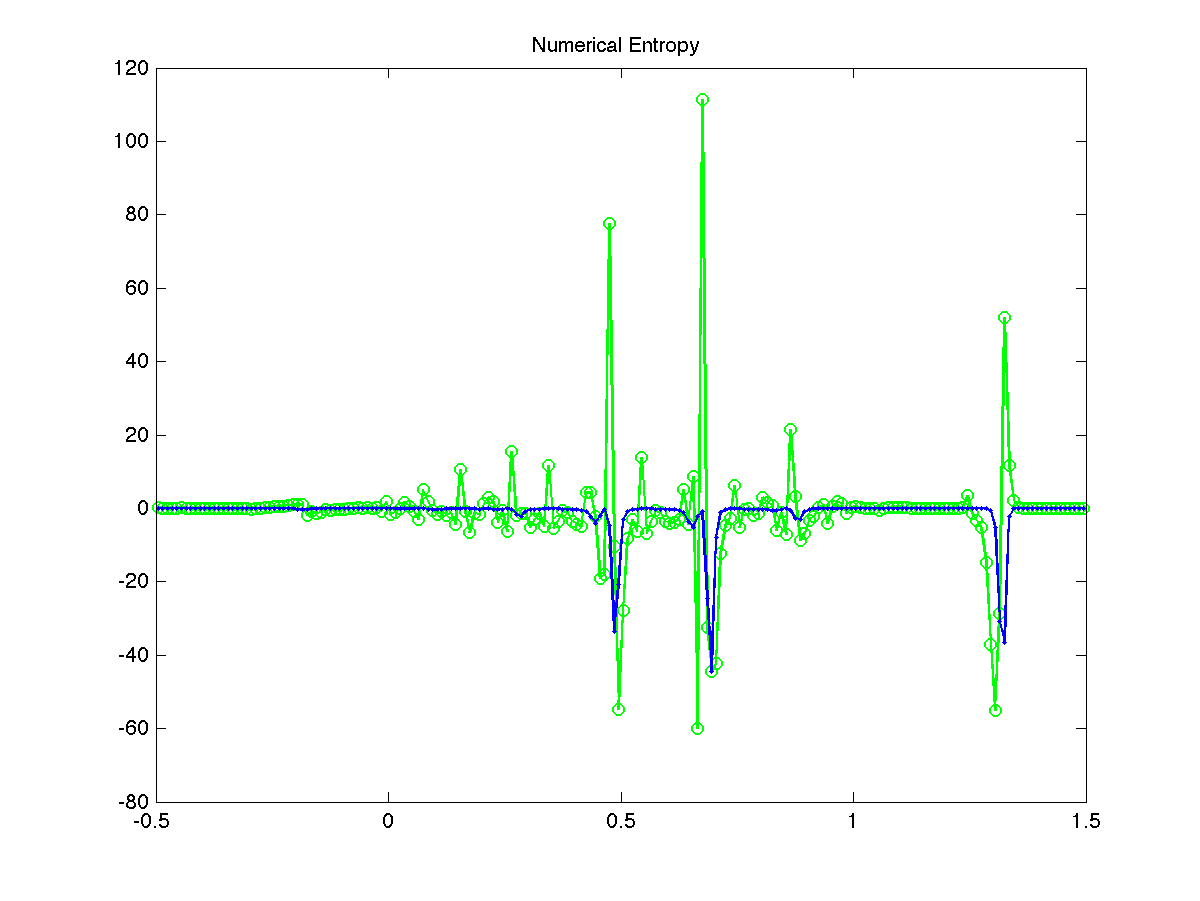}
\\
Third order
&
Fourth order
\\
\includegraphics[width=0.45\linewidth]{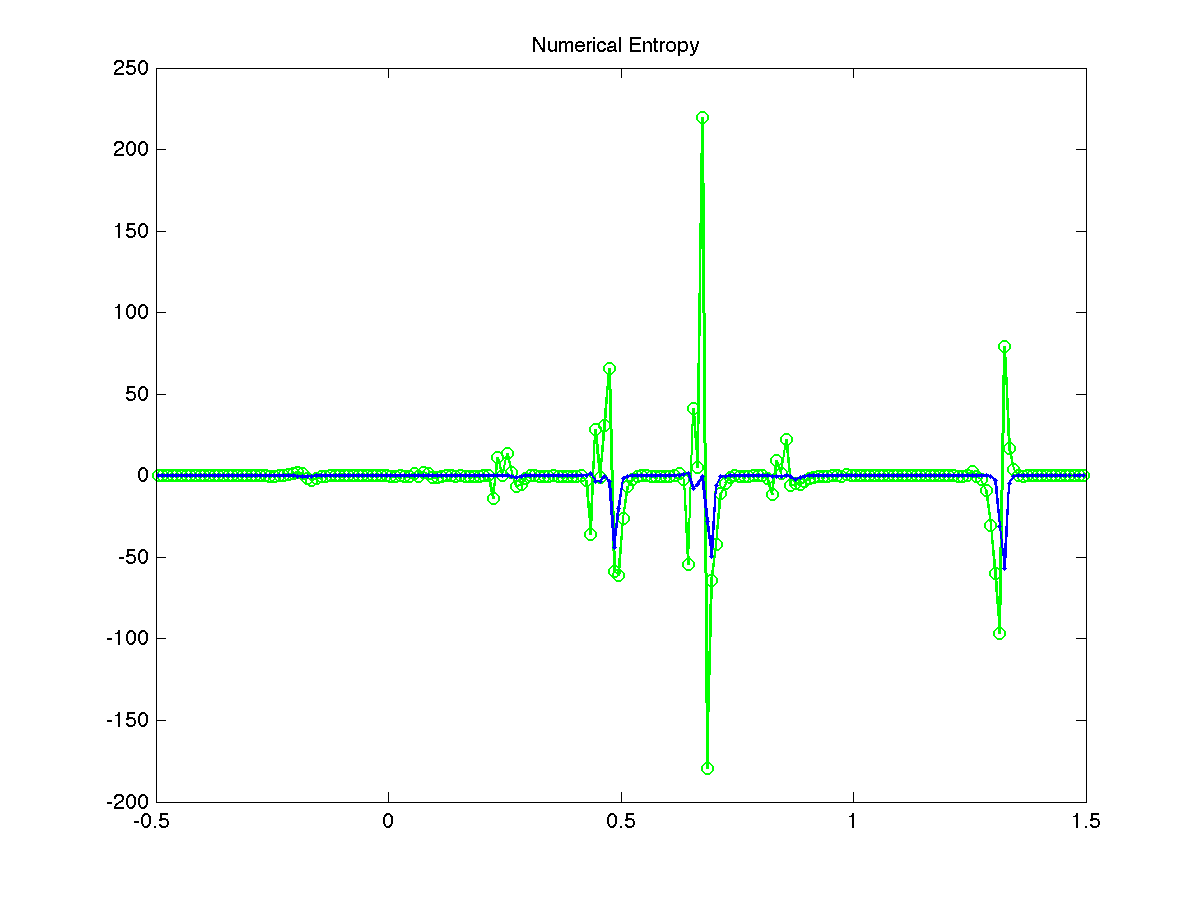}
&
\includegraphics[width=0.45\linewidth]{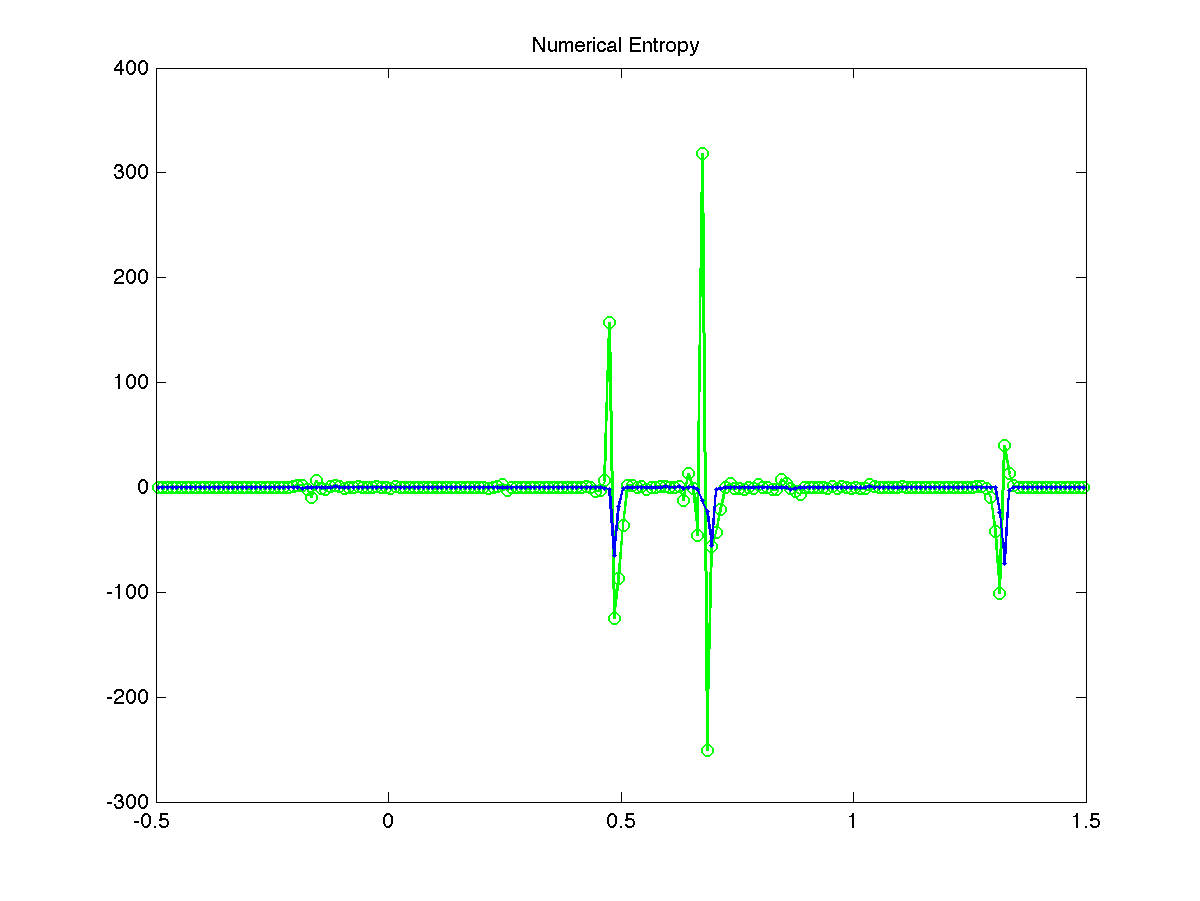}
\end{tabular}
\end{center}
\caption{Comparison of the numerical entropy production with two different numerical entropy fluxes.}
\label{fig:nument:compareflux}
\end{figure}

\medskip
Finally, we wish to illustrate the importance of choosing the numerical entropy flux customized on the numerical flux used by the scheme, as in \eqref{eq:numentflux}.
Figure \ref{fig:nument:compareflux} shows the numerical entropy production on the test \eqref{eq:SinCanal} computed with the numerical entropy flux of \eqref{eq:numentflux} (green line with circles) and with the numerical entropy flux $\Psi(U^-,U^+)=\tfrac12 (U^-+U^+)$ (blue line with dots). Note that also the alternative flux considered here is consistent with the exact entropy flux $\psi$ and therefore will provide entropy residuals with the same rate of decay of the local error of the scheme.

However, in all cases, it is clear that using the local Lax-Friedrichs flux for both the conservation law and the computation of the numerical entropy flux leads to much smaller positive overshoots in the numerical entropy production and thus a much more reliable error indicator.

\section{Conclusions}
In this work we have derived formulas for high order schemes for balance laws on non-uniform grids. It includes the extension of the third order compact WENO reconstruction of \cite{LPR01} to non uniform grids and high order reconstructions to compute the cell average of the source term, needed by high order finite volume schemes on balance laws.
Farther, we illustrate how well balancing on equilibrium solutions can be enforced for high order schemes on irregular grids. 

We also include the extension of the entropy indicator we proposed in \cite{PS:entropy} and \cite{P:entropy} to the case of balance laws. The proofs given in \cite{PS:entropy} carry over to the case of balance laws with geometric source terms, and prove that the entropy indicator provides a measure of the local truncation error on smooth flows, and it reliably selects the location of discontinuities.

Several numerical tests are included, to show the achievement of the expected accuracy of the schemes proposed, even on extremely irregular grids, and the improvement obtained with ad-hoc chosen grids. 

Future work on this topic will be dedicated to the construction of adaptive cartesian grids of octree type, driven by the entropy error indicator, for balance laws, with particular attention on the enforcement of equilibrium solutions at the discrete level.

\begin{acknowledgements}
This work was supported by 
\textquotedblleft National Group for Scientific Computation (GNCS-INDAM) 
\textquotedblright
\end{acknowledgements}

\bibliographystyle{spmpsci}
\bibliography{shentropy}
\nocite{*}

\end{document}